# GRIDAPROMS.JL: EFFICIENT REDUCED ORDER MODELLING IN THE JULIA PROGRAMMING LANGUAGE

NICHOLAS MUELLER† AND SANTIAGO BADIA†

ABSTRACT. In this paper, we introduce GridapROMs, a Julia-based library for the numerical approximation of parameterized partial differential equations (PDEs) using a comprehensive suite of linear reduced order models (ROMs). The library is designed to be extendable and productive, leveraging an expressive high-level API built on the Gridap PDE solver backend, while achieving high performance through Julia's just-in-time compiler and advanced lazy evaluation techniques. GridapROMs is PDE-agnostic, enabling its application to a wide range of problems, including linear, nonlinear, single-field, multi-field, steady, and unsteady equations. This work details the library's key innovations, implementation principles, and core components, providing usage examples and demonstrating its capabilities by solving a fluid dynamics problem modeled by the Navier-Stokes equations in a 3D geometry.

**Program summary**
*Program Title:* GridapROMs.jl (version 1.0)
*Developer's repository link:* https://github.com/Gridap/GridapROMs.jl
*Licensing provisions:* MIT license
*Programming language:* Julia
*Nature of problem:* Numerical simulation of parameterized PDEs, including linear, nonlinear, single-field, multi-field, steady, and unsteady problems. Classical full-order models are computationally expensive, requiring intensive computations for each parameter configuration.
*Solution method:* GridapROMs approximates the parameter-to-solution map using linear reduced order models. It constructs a reduced basis from the tangent hyperplane to the solution manifold and applies a (Petrov-)Galerkin projection to the full-order equations. Nonaffine parameter dependencies in the residual and/or Jacobian are efficiently handled using hyper-reduction techniques.

## 1. INTRODUCTION

Conventional high-fidelity (HF) solvers for parametric partial differential equations (PDEs) employ fine discretizations for the numerical integration of weak formulations. This process results in the assembly of large systems of equations – referred to as full-order models (FOMs) – which are then solved using appropriate numerical schemes. Even with the support of parallel toolboxes, such solvers incur substantial computational costs [1], especially in multi-query settings where solutions must be computed for many parameter values. These costs are further exacerbated in unsteady problems. To address this challenge, reduced-order models (ROMs) have emerged as efficient alternatives, offering low-dimensional approximation spaces for solving PDEs at reduced computational cost. Among these, the reduced basis (RB) method stands out as a prominent data-driven, projection-based ROM, leveraging HF solution snapshots to compute reduced subspaces and applying Galerkin projection to the FOM equations. For linearly reducible problems [2, 3], RB methods demonstrate remarkable accuracy at a fraction of the computational cost, particularly in unsteady scenarios [4, 5]. Despite their popularity, open-source implementations of RB methods remain scarce. Notable examples include redBKit [2], the RB module in libMesh [6], RBmatlab, Dune-RB [7], libROM [8], RBniCS [9] and pyMOR [10], have also been developed, with pyMOR being particularly well-known. We also mention [11], which presents a simple Python implementation for time-dependent ROMs. These libraries often rely on computationally-intensive kernels implemented in precompiled languages like C++, such as PDE solvers like DOLFIN/FEniCS [12] and Deal.II [13] or NumPY for linear algebra, while providing high-level scripting interfaces in interpreted languages like Python for ease of use. Although this approach ensures computational efficiency, it introduces the two-language barrier, requiring users to navigate between a high-level scripting language for usability and a low-level language for performance-critical tasks. This duality can complicate development workflows, hinder extendibility, and reduce overall productivity.

GridapROMs is a novel, open-source RB library written entirely in the Julia programming language, designed to overcome the limitations of existing libraries, particularly the two-language barrier. Julia seamlessly combines the user-friendliness of interpreted languages like Python with the high performance of compiled languages like C++, leveraging its just-in-time (JIT) compiler to produce highly optimized native machine code tailored to runtime data types [14]. GridapROMs relies on the PDE solver Gridap [15, 16] for the HF computations, which is also entirely written in Julia.

†SCHOOL OF MATHEMATICS, MONASH UNIVERSITY, CLAYTON, VICTORIA 3800, AUSTRALIA
*E-mail addresses*: nicholas.mueller@monash.edu, santiago.badia@monash.edu.

Gridap combines high performance with a user-friendly interface that closely resembles mathematical notation. Moreover, the Julia package manager provides a robust ecosystem of interoperable libraries, streamlining integration and enhancing productivity in scientific computing (e.g., JuMP [17] for mathematical optimization or DifferentialEquations [18] for solving ordinary differential equations).

Building upon the expressive API of Gridap, GridapROMs enables a seamless definition of the FOM and its reduced counterpart. The library achieves high performance by leveraging advanced programming techniques, and state-of-the-art algorithms for reduced order modeling. A distinctive feature of GridapROMs is its use of *lazy evaluations* for HF quantities, taking advantage of the capabilities provided by Gridap. Lazy evaluations enable the composition of functions without eagerly computing intermediate results, significantly enhancing efficiency. As detailed in Subsection 2.3, this approach is crucial for alleviating the computational cost of collecting HF snapshots, a task commonly dubbed as a "computational bottleneck" within the ROM community. Additionally, GridapROMs implements a distributed-in-memory framework to efficiently generate solution, residual, and Jacobian snapshots for parameterized PDEs, enabling scalable reduced-order modeling in parallel computing environments. After collecting snapshots, the reduced subspace is computed using an appropriate rank-reduction technique. GridapROMs supports a variety of compression strategies, including standard truncated proper orthogonal decomposition (TPOD) methods [2, 19] and efficient randomized algorithms [20, 21]. The library also allows the construction of local approximation subspaces – rather than standard global ones – which have proven highly effective for compressing poorly reducible problems [22–26]. Furthermore, GridapROMs incorporates advanced tensor train (TT)-based ROM techniques [27], particularly suited for high-dimensional problems, representing a significant advancement over standard reduced-order modeling libraries. Subsequently, we define a Galerkin projection operator to map HF quantities onto the reduced subspace. To further enhance efficiency, we employ state-of-the-art hyper-reduction techniques, such as those in [3, 5, 28], which significantly reduce the computational complexity of the FOM. These steps comprise the *offline phase* of the method, which is computationally intensive but performed only once. In the subsequent *online phase*, we efficiently compute the corresponding RB solution for any new parameter, utilizing precomputed reduced operators and hyper-reduction techniques.

The paper is organized as follows. Section 2 introduces the mathematical foundations of the RB method, followed by a detailed discussion of the design principles, core components, and a usage example of GridapROMs. Section 3 highlights some of the library's more advanced features with illustrative code snippets, including stabilization techniques for ROMs applied to saddle-point problems, TT-based ROM methods, local subspaces, and the generation of snapshots using distributed-in-memory parallelism. Section 4 demonstrates the application of GridapROMs to a fluid dynamics problem modeled by a transient Navier-Stokes equation in a 3D geometry, with parameterizations affecting both the Reynolds number and boundary conditions. Section 5 provides concluding remarks and outlines potential directions for future development. The paper also includes two appendices: Appendix A presents the code used to run the application of Section 4, while Appendix B illustrates how GridapROMs can be leveraged to solve a problem on a parameterized geometry using a TT-based local subspace strategy.

## 2. Formulation and implementation details

We begin this section with the mathematical formulation of the RB method for parameterized PDEs, first for steady-state problems, then for time-dependent ones. Subsequently, we describe the design principles and the main abstractions in GridapROMs. We conclude the section by providing a usage example of the library.

### 2.1. Mathematical formulation.

Consider a parameterized PDE defined on a domain $\Omega^\mu \subset \mathbb{R}^d$, where $d = 2, 3$, and characterized by a parameter $\boldsymbol{\mu} \in \mathcal{D} \subset \mathbb{R}^p$, with $\mathcal{D}$ representing the parameter space. The general form of such a PDE is given by

$$\text{find } u^\mu = u^\mu(\underline{x}) \in \mathcal{U} \text{ such that } \mathcal{A}^\mu(u^\mu) = 0, \quad \underline{x} \in \Omega^\mu, \tag{1}$$

subject to a set of boundary conditions on $\partial\Omega^\mu$. Here, $\mathcal{U}$ is a space of sufficiently smooth functions defined on $\Omega^\mu$, $u^\mu$ is the unknown solution, and $\mathcal{A}^\mu$ is a (nonlinear) differential operator. The superindex $\mu$ indicates the dependence on the parameter $\boldsymbol{\mu}$. Note that the domain $\Omega^\mu$ may vary with $\boldsymbol{\mu}$, allowing for shape parameters to influence the geometry of the problem, making this the most general form of a parametric PDE. For simplicity, we assume $\Omega$ to be fixed henceforth. Eq. (1) is commonly referred to as the strong form of the PDE. To proceed with the finite element (FE) discretization, we introduce a quasi-uniform partition of $\Omega$, denoted as $\mathcal{T}_h$, where $h$ represents the mesh size, and define a pair of trial and test FE spaces $(\mathcal{U}_h^\mu, \mathcal{V}_h)$ on $\mathcal{T}_h$. Note that the trial space is generally characterized by a parametric dependence via the Dirichlet boundary conditions. Let $v_h \in \mathcal{V}_h$ be an arbitrary test function and $u_h^\mu \in \mathcal{U}_h^\mu$ the FE approximation of the unknown. The weak formulation corresponding to (1) is given by

$$\text{find } u_h^\mu = u_h^\mu(\underline{x}) \in \mathcal{U}_h^\mu \text{ such that } a_h^\mu(u_h^\mu, v_h) = 0, \quad \forall v_h \in \mathcal{V}_h, \quad \underline{x} \in \Omega. \tag{2}$$

The nonlinear form $a_h^\mu$ is derived by multiplying $\mathcal{A}^\mu$ by $v_h$ and integrating by parts over $\Omega$. For the well-posedness of (2), it is sufficient to assume that the operator $\mathcal{A}^\mu$ is continuously differentiable with respect to $u_h^\mu$ (ensuring the existence of a solution) and that the problem satisfies a small data assumption (ensuring uniqueness). Hereafter, we operate under the assumption of well-posedness. Leveraging the differentiability of $\mathcal{A}^\mu$, we linearize (2) and solve the resulting problem iteratively using the Newton-Raphson method. The linearized problem can be expressed algebraically as

$$\text{given } \boldsymbol{w}_h^{(0)} \in \mathbb{R}^N, \text{ compute } \boldsymbol{J}_h^\mu(\boldsymbol{w}_h^{(k)})\delta\boldsymbol{w}_h^{(k)} = -\boldsymbol{r}_h^\mu(\boldsymbol{w}_h^{(k)}), \text{ and update } \boldsymbol{w}_h^{(k+1)} = \boldsymbol{w}_h^{(k)} + \delta\boldsymbol{w}_h^{(k)}, \text{ for } k = 1, 2, \ldots \tag{3}$$

When a stopping criterion is met, for example

$$\|\delta \boldsymbol{w}_h^{(k)}\| < \varepsilon,$$

we then set $\boldsymbol{u}_h^{\mu} = \boldsymbol{w}_h^{(k+1)}$, where $\varepsilon$ represents a sufficiently small tolerance. In (3), $\boldsymbol{J}_h^{\mu}$ is the $N \times N$ nonlinear Jacobian matrix, obtained from the numerical integration of the Fréchet derivative [29, 30] of $a_h^{\mu}$, while $\boldsymbol{r}_h^{\mu}$ is the $N$-dimensional residual vector, derived from the numerical integration of $a_h^{\mu}$. Here, $N$ denotes the number of full-order degrees of freedom (DOFs) in the problem. Since $N$ is typically very large in practical HF applications, ROMs are designed to replace the FOM system (3) with a significantly smaller system of equations that still accurately approximate the solution.

The RB method is one of the most widely recognized projection-based ROMs. It begins by solving the FOM for a set of offline realizations $\boldsymbol{\mu}_{\texttt{off}} \subset \mathcal{D}$ to generate a tensor of snapshots. From these snapshots, a reduced basis is constructed using a low-rank approximation algorithm, such as the standard TPOD for steady problems, a space-time TPOD for transient problems [4, 5, 31], or a tensor decomposition like tensor train SVD (TT-SVD) [32, 33]. Letting $\boldsymbol{\Phi} \in \mathbb{R}^{N \times n}$ represent the reduced basis obtained from one of these techniques, the RB approximation is expressed as

$$\boldsymbol{u}_h^{\mu} \approx \boldsymbol{u}_n^{\mu} = \boldsymbol{\Phi} \widehat{\boldsymbol{u}}^{\mu}, \tag{4}$$

where $\widehat{\boldsymbol{u}}^{\mu}$ is the $n$-dimensional vector of unknown coefficients in the reduced basis, with $n \ll N$. We then define the reduced trial and test spaces $(\mathcal{U}_n, \mathcal{V}_n)$, where $\mathcal{U}_n = \text{Col}(\boldsymbol{\Phi}) \subset \mathcal{U}_h^{\mu}$ and $\mathcal{V}_n = \text{Col}(\boldsymbol{\Psi}) \subset \mathcal{V}_h$, to derive a reduced version of (2). Here, Col denotes the column space of a matrix, and $\boldsymbol{\Psi} \in \mathbb{R}^{N \times n}$ is a (full-column rank) matrix whose expression will be defined later. Note that $\mathcal{U}_n$ does not feature a $\boldsymbol{\mu}$-dependence, unlike its full-order counterpart $\mathcal{U}_h^{\mu}$. Indeed, $\boldsymbol{\Phi}$ is usually computed from the free values of the solution snapshots, and the Dirichlet datum is simply added to the ROM approximation as a lifting term. Since $\mathcal{U}_n \subset \mathcal{U}_h^{\mu}$, each RB vector $\Phi_i$ can be expressed in terms of the FE basis $\{\varphi_j\}_{j=1}^N$ as

$$\Phi_i = \sum_{j=1}^{N} \varphi_j \Phi_{j,i}.$$

Now let us refer to an arbitrary reduced test function as $v_n \in \mathcal{V}_n$, and to $u_n^{\mu} \in \mathcal{U}_n$ as the FE function

$$u_n^{\mu}(\underline{x}) = \sum_{i=1}^{n} \Phi_i(\underline{x}) \widehat{u}_i^{\mu} = \sum_{i=1}^{n} \left( \sum_{j=1}^{N} \varphi_j(\underline{x}) \Phi_{j,i} \right) \widehat{u}_i^{\mu}. \tag{5}$$

If we require the approximant $u_n^{\mu}$ to satisfy the weak formulation (2) for any $v_n$, we get the Petrov-Galerkin projection equation:

$$\text{find } u_n^{\mu} = u_n^{\mu}(\underline{x}) = \sum_{i=1}^{n} \Phi_i(\underline{x}) \widehat{u}_i^{\mu} \in \mathcal{U}_n \text{ such that } a_h^{\mu}(u_n^{\mu}, v_n) = 0, \quad \forall v_n \in \mathcal{V}_n, \quad \underline{x} \in \Omega. \tag{6}$$

We can algebraically write the expression above as

$$\text{given } \widehat{\boldsymbol{w}}^{(0)} \in \mathbb{R}^n, \text{ compute } \widehat{\boldsymbol{J}}^{\mu}(\widehat{\boldsymbol{w}}^{(k)}) \delta \widehat{\boldsymbol{w}}^{(k)} = -\widehat{\boldsymbol{r}}^{\mu}(\widehat{\boldsymbol{w}}^{(k)}), \text{ and update } \widehat{\boldsymbol{w}}^{(k+1)} = \widehat{\boldsymbol{w}}^{(k)} + \delta \widehat{\boldsymbol{w}}^{(k)} \tag{7}$$

where

$$\widehat{\boldsymbol{J}}^{\mu}(\widehat{\boldsymbol{w}}) = \boldsymbol{\Psi}^T \boldsymbol{J}_h^{\mu}(\widehat{\boldsymbol{w}}) \boldsymbol{\Phi}; \qquad \widehat{\boldsymbol{r}}^{\mu}(\widehat{\boldsymbol{w}}) = \boldsymbol{\Psi}^T \boldsymbol{r}_h^{\mu}(\widehat{\boldsymbol{w}}). \tag{8}$$

While the framework readily supports Petrov-Galerkin projections, we assume from now that $\mathcal{V}_n \equiv \mathcal{U}_n$ for simplicity in the exposition. It is worth noting that Petrov-Galerkin formulations offer advantages over Galerkin projections only in specific scenarios (e.g., [4, 19]).

A key principle of the RB method (and ROMs in general) is the separation of computations into an offline phase and an online phase. During the offline phase, we construct $\boldsymbol{\Phi}$ and precompute the projected quantities in (8). While these operations are computationally intensive, they are performed only once. In the online phase, we efficiently compute the reduced solution $\widehat{\boldsymbol{u}}^{\mu}$ by solving (7) for any given $\boldsymbol{\mu}$, with a cost independent of $N$. However, in many practical applications, the full-order left-hand side (LHS) and right-hand side (RHS) depend on $\boldsymbol{\mu}$, making it infeasible to fully precompute (8) offline. To address this, hyper-reduction techniques are employed to approximate $\boldsymbol{J}_h^{\mu}$ and $\boldsymbol{r}_h^{\mu}$ using affine decompositions:

$$\boldsymbol{J}_h^{\mu}(\widehat{\boldsymbol{w}}) \approx \boldsymbol{J}_{n,n}^{\mu}(\widehat{\boldsymbol{w}}) = \sum_{i=1}^{n^{\boldsymbol{J}}} \boldsymbol{\Phi}_i^{\boldsymbol{J}} \widehat{J}_i^{\mu}(\widehat{\boldsymbol{w}}); \qquad \boldsymbol{r}_h^{\mu}(\widehat{\boldsymbol{w}}) \approx \boldsymbol{r}_n^{\mu}(\widehat{\boldsymbol{w}}) = \sum_{i=1}^{n^{\boldsymbol{r}}} \boldsymbol{\Phi}_i^{\boldsymbol{r}} \widehat{r}_i^{\mu}(\widehat{\boldsymbol{w}}). \tag{9}$$

Here, $\boldsymbol{\Phi}_i^{\boldsymbol{J}} \in \mathbb{R}^{N \times N}$ for $i = 1, \ldots, n^{\boldsymbol{J}}$ and $\boldsymbol{\Phi}_i^{\boldsymbol{r}} \in \mathbb{R}^N$ for $i = 1, \ldots, n^{\boldsymbol{r}}$ denote the bases for the $n^{\boldsymbol{J}}$- and $n^{\boldsymbol{r}}$-dimensional subspaces approximating the manifolds of parameterized Jacobians and residuals, respectively. Additionally, $\widehat{\boldsymbol{J}}^{\mu} \in \mathbb{R}^{n^{\boldsymbol{J}}}$ and $\widehat{\boldsymbol{r}}^{\mu} \in \mathbb{R}^{n^{\boldsymbol{r}}}$ represent the reduced coefficients of $\boldsymbol{J}_h^{\mu}$ and $\boldsymbol{r}_h^{\mu}$ with respect to their corresponding bases. Common hyper-reduction techniques include discrete empirical interpolation method (DEIM) [34], its matrix version (MDEIM) [3, 5], the subspace-orthogonal point selection (S-OPT) [28], and other collocation methods described in [35]. In the presentation of the methodology below, we consider a MDEIM-based hyper-reduction, which we briefly outline with the aid of Fig. 1:

(1) Collect snapshots $\{\boldsymbol{J}_h^{\mu}\}_{\boldsymbol{\mu} \in \boldsymbol{\mu}_{\texttt{off}}}$ and horizontally concatenate their vectors of nonzero entries. This step assumes that all Jacobian snapshots share the same sparsity pattern, enabling the concatenation.

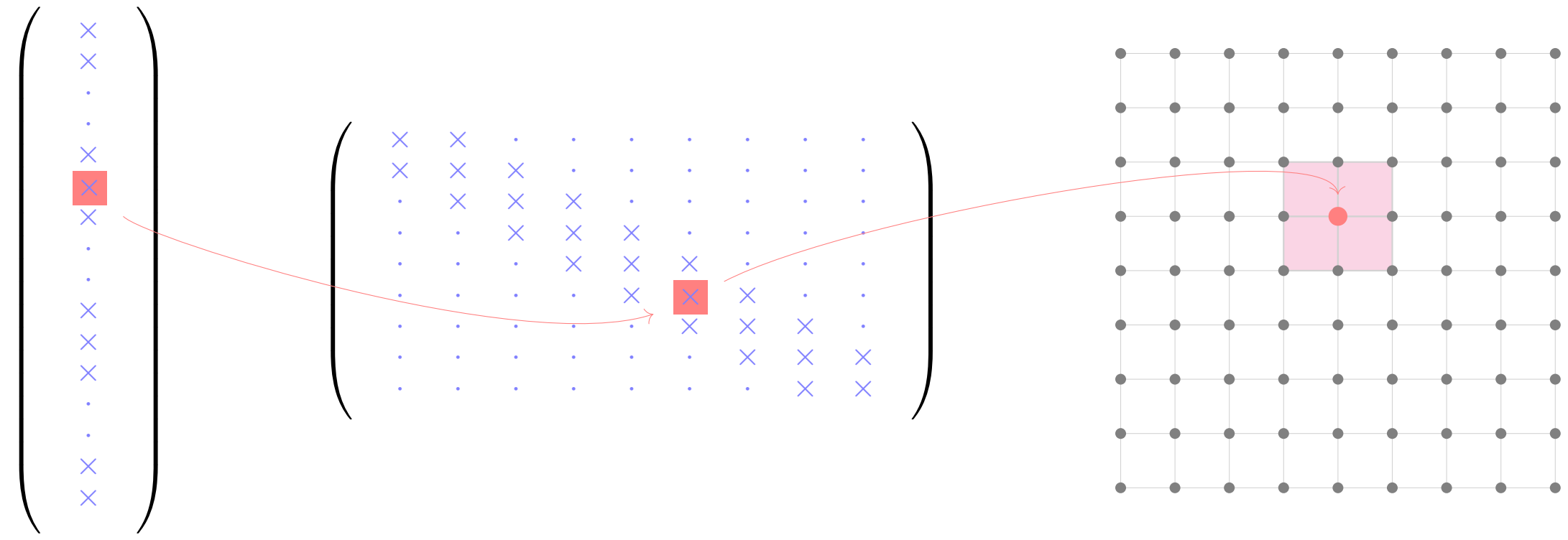

FIGURE 1. Graphical representation of the *reduced integration domain* in the MDEIM approximation of the Jacobian. The middle figure illustrates $\boldsymbol{\Phi}_i^{\boldsymbol{J}}$, the $i$th component of the Jacobian basis, which retains the sparsity pattern of $\boldsymbol{J}_h^{\mu}$. On the left, the corresponding vector of nonzero entries is shown, while on the right, the FE mesh of the problem is depicted. The first arrow represents a bijective mapping between a nonzero entry of $\boldsymbol{\Phi}_i^{\boldsymbol{J}}$ and a row-column index pair. The second arrow establishes another bijective mapping that associates each row-column index pair with a set of integration cells. The MDEIM procedure applied to the Jacobian basis identifies a list of FE cells, defining the reduced integration domain.

(2) Apply a rank-reduction technique (e.g., TPOD) to the concatenated snapshots to extract the basis $\boldsymbol{\Phi}_z^{\boldsymbol{J}} \in \mathbb{R}^{N_z \times n^{\boldsymbol{J}}}$, where $n^{\boldsymbol{J}} \ll N_z$ and $N_z$ is the number of nonzero entries.
(3) Construct a vector of interpolation indices $\mathcal{J} = [j_1, \dots, j_{n^{\mathcal{J}}}] \subset \{1, \dots, N_z\}^{n^{\boldsymbol{J}}}$ iteratively. For each column $\boldsymbol{\Phi}_z^{\boldsymbol{J}}[:, i]$ of the basis, select the index $j_i$ corresponding to the row that maximizes a residual-like estimator, as described in [34]. The selected indices are marked in red in Fig. 1.
(4) During the online phase, compute the reduced coefficients $\widehat{\boldsymbol{J}}^{\mu}$ for any parameter $\boldsymbol{\mu}$ using the formula:

$$\widehat{\boldsymbol{J}}^{\mu} = \boldsymbol{\Phi}_z^{\boldsymbol{J}}[\mathcal{J}, :]^{-1} \boldsymbol{J}_z^{\mu}[\mathcal{J}], \tag{10}$$

where $\boldsymbol{J}_z^{\mu}$ is the vector of nonzero entries of $\boldsymbol{J}_h^{\mu}$. The term $\boldsymbol{J}_z^{\mu}[\mathcal{J}]$ is efficiently computed by restricting the cell-wise integration and assembly routines to the FE cells identified by $\mathcal{J}$, as illustrated in Fig. 1.

After determining the affine decompositions in (9), we solve the hyper-reduced ROM system by substituting these terms into (7):

$$\text{given } \widehat{\boldsymbol{w}}^{(0)} \in \mathbb{R}^n, \text{ compute } \bar{\boldsymbol{J}}(\widehat{\boldsymbol{w}}^{(k)})\delta\widehat{\boldsymbol{w}}^{(k)} = -\bar{\boldsymbol{r}}^{\mu}(\widehat{\boldsymbol{w}}^{(k)}), \text{ and update } \widehat{\boldsymbol{w}}^{(k+1)} = \widehat{\boldsymbol{w}}^{(k)} + \delta\widehat{\boldsymbol{w}}^{(k)} \tag{11}$$

where

$$\bar{\boldsymbol{J}}^{\mu}(\widehat{\boldsymbol{w}}) = \sum_{i=1}^{n^{\boldsymbol{J}}} \boldsymbol{\Phi}^T \boldsymbol{\Phi}_i^{\boldsymbol{J}} \boldsymbol{\Phi} \widehat{J}_i^{\mu}(\widehat{\boldsymbol{w}}); \qquad \bar{\boldsymbol{r}}^{\mu}(\widehat{\boldsymbol{w}}) = \sum_{i=1}^{n^{\boldsymbol{r}}} \boldsymbol{\Phi}^T \boldsymbol{\Phi}_i^{\boldsymbol{r}} \widehat{r}_i^{\mu}(\widehat{\boldsymbol{w}}). \tag{12}$$

The computation of $\{\boldsymbol{\Phi}^T \boldsymbol{\Phi}_i^{\boldsymbol{J}} \boldsymbol{\Phi}\}_{i=1}^{n^{\boldsymbol{J}}}$ and $\{\boldsymbol{\Phi}^T \boldsymbol{\Phi}_i^{\boldsymbol{r}}\}_{i=1}^{n^{\boldsymbol{r}}}$ constitutes the bulk of the offline operations required to solve (11). Although these computations are resource-intensive, they are performed only once during the offline phase, as they are independent of $\boldsymbol{\mu}$. During the online phase, the only $\boldsymbol{\mu}$-dependent terms to compute are the reduced coefficients $\widehat{\boldsymbol{J}}^{\mu}$ and $\widehat{\boldsymbol{r}}^{\mu}$. Using the procedure outlined in Fig. 1, these coefficients can be computed efficiently at a cost independent of $N$. Once obtained, the terms in (12) are assembled, and the Newton-Raphson iterations in (11) are solved. This step is computationally inexpensive, as it involves inverting a matrix of size $n \times n$ at each iteration.

2.2. **Mathematical formulation of time-dependent problems.** This subsection introduces the benchmark ROM for time-dependent, nonlinear, parameterized PDEs. We start by presenting the weak formulation: given an initial condition

$$u_h^{\mu}(\underline{x}, 0) = u_0^{\mu}(\underline{x}) \quad \underline{x} \in \Omega,$$

find $u_h^{\mu} = u_h^{\mu}(\underline{x}, t) \in \mathcal{U}_h^{\mu}$ such that

$$\left(\frac{\partial u_h^{\mu}}{\partial t}, v_h\right) + a_h^{\mu}(u_h^{\mu}, v_h) = 0, \quad \forall v_h \in \mathcal{V}_h, \quad (\underline{x}, t) \in \Omega \times (0, T],$$

and subject to appropriate boundary conditions on $\partial\Omega \times (0, T]$. The solution is defined over the space-time domain $\Omega \times [0, T]$, with $T > 0$. To derive the space-time FOM, we discretize the temporal domain into a uniform partition $\{t_n\}_{n=0}^{N_t}$, where $t_n = n\Delta t$ and $\Delta t = T/N_t$ denotes the time-step size. A time marching scheme is then applied to compute the fully discrete solution. For instance, the Backward Euler (BE) method at the $k$th iteration is expressed as:

$$\text{given } \boldsymbol{w}_{(n)}^{(0)} \in \mathbb{R}^N, \text{ compute } \Delta t^{-1}\boldsymbol{M}(\delta\boldsymbol{w}_{(n)}^{(k)} - \boldsymbol{u}_{(n-1)}^{\mu}) + \boldsymbol{J}_h^{\mu}(\boldsymbol{w}_{(n)}^{(k)})\delta\boldsymbol{w}_{(n)}^{(k)} = -\boldsymbol{r}_h^{\mu}(\boldsymbol{w}_{(n)}^{(k)}), \text{ update } \boldsymbol{w}_{(n)}^{(k+1)} = \boldsymbol{w}_{(n)}^{(k)} + \delta\boldsymbol{w}_{(n)}^{(k)}. \tag{13}$$

Here, the variable $\boldsymbol{u}^{\mu}_{(n)}$ represents the FOM solution at the $n$th time step. By definition of the initial condition, we have $\boldsymbol{u}^{\mu}_{(0)} = \boldsymbol{u}^{\mu}_{0}$, where $\boldsymbol{u}^{\mu}_{0}$ corresponds to the nodal values of the initial condition. Equation (13) can be reformulated as the following tridiagonal block system:

$$
\begin{bmatrix}
\Delta t^{-1}\boldsymbol{M} + \boldsymbol{J}^{\mu}_{h}(\boldsymbol{w}^{(k)}_{(1)}) & \mathbf{0} & & \\
-\Delta t^{-1}\boldsymbol{M} & \Delta t^{-1}\boldsymbol{M} + \boldsymbol{J}^{\mu}_{h}(\boldsymbol{w}^{(k)}_{(2)}) & \ddots & \\
 & \ddots & \ddots & \mathbf{0} \\
 & & -\Delta t^{-1}\boldsymbol{M} & \Delta t^{-1}\boldsymbol{M} + \boldsymbol{J}^{\mu}_{h}(\boldsymbol{w}^{(k)}_{(N_t)})
\end{bmatrix}
\begin{bmatrix}
\delta\boldsymbol{w}^{(k)}_{1} \\ \delta\boldsymbol{w}^{(k)}_{(2)} \\ \vdots \\ \delta\boldsymbol{w}^{(k)}_{(N_t)}
\end{bmatrix}
= -
\begin{bmatrix}
\Delta t^{-1}\boldsymbol{M}\boldsymbol{w}^{(k)}_{(0)} + \boldsymbol{r}^{\mu}_{h}(\boldsymbol{w}^{(k)}_{(1)}) \\ \boldsymbol{r}^{\mu}_{h}(\boldsymbol{w}^{(k)}_{(2)}) \\ \vdots \\ \boldsymbol{r}^{\mu}_{h}(\boldsymbol{w}^{(k)}_{(N_t)})
\end{bmatrix}
\tag{14}
$$

We refer to (14) as the FOM for transient applications. While the system is not explicitly assembled in practice, expressing it in this form is insightful, as it highlights that by introducing a space-time variable

$$\boldsymbol{w}_{h\Delta} = [\boldsymbol{w}_{(1)}, \ldots, \boldsymbol{w}^{T}_{(N_t)}] \in \mathbb{R}^{N \cdot N_t}$$

we can compactly rewrite the transient FOM as

$$\boldsymbol{J}^{\mu}_{\Delta}(\boldsymbol{w}^{(k)}_{h\Delta})\delta\boldsymbol{w}^{(k)}_{h\Delta} = -\boldsymbol{r}^{\mu}_{\Delta}(\boldsymbol{w}^{(k)}_{h\Delta}), \quad \text{where} \quad \boldsymbol{J}^{\mu}_{\Delta} \in \mathbb{R}^{N \cdot N_t \times N \cdot N_t},\ \boldsymbol{r}^{\mu}_{\Delta} \in \mathbb{R}^{N \cdot N_t}.$$

Now, we consider a space-time projection operator $\boldsymbol{\Phi} \in \mathbb{R}^{N \cdot N_t \times n}$, which can be built by employing either the space-time reduced basis (ST-RB) method proposed in [4, 5, 36], or the tensor train reduced basis (TT-RB) procedure in [27]. By following the procedure outlined in Eqs. (5)-(6), we can write a transient ROM that reads exactly as Eq. (7). In practice, the transient ROM eliminates the time marching. For the approximation of the space-time Jacobians and residuals, we can employ a space-time hyper-reduction introduced in [5, 27]. In essence, we consider the space-time bases for the Jacobians and residuals

$$\boldsymbol{\Phi}^{\boldsymbol{J}}_{i} \in \mathbb{R}^{N \cdot N_t \times N \cdot N_t}\ \forall\, i = 1, \ldots, n^{\boldsymbol{J}}; \qquad \boldsymbol{\Phi}^{\boldsymbol{r}}_{i} \in \mathbb{R}^{N \cdot N_t}\ \forall\, i = 1, \ldots, n^{\boldsymbol{r}} \tag{15}$$

and substituting the resulting affine decompositions (see (9)) into the transient ROM leads to the same hyper-reduced system (11). The structures in (15) can be efficiently computed by following the methodology outlined in [4, 5]. Consequently, instead of solving a FOM that involves a Newton-Raphson loop nested within a time marching scheme, the ROM reduces the problem to efficient space-time Newton-Raphson iterations. As demonstrated in the numerical results in Section 4, this approach offers significant computational speedup for transient applications.

2.3. **Implementation principles.** A performant RB library must ensure efficiency in both offline and online operations. While achieving high performance in the online phase is relatively straightforward, the offline phase often represents a significant computational challenge, commonly referred to in the ROM community as the "computational bottleneck." To address this, it is essential to focus on:

- An efficient generation and storage of the snapshots.
- State-of-the-art reduction algorithms for the computation of the RBs.

In this subsection, we focus on the first point, which is significantly more challenging to achieve. For now, we assume the computations are performed serially, whereas we address the parallel case in Subsection 3.4. GridapROMs extends Gridap [15, 16] with an efficient mechanism to perform FE tasks (e.g., integration, assembly, and solve) for any number of parameters. The key to this efficiency lies in the use of lazy operations on parametric HF quantities. Lazy operations, a well-established concept in functional programming, differ from standard *eager* operations by deferring both allocation and computation of output entries until explicitly required. Instead of allocating an output structure and populating its entries immediately, lazy operations return a wrapper structure encapsulating the operator and its arguments. This wrapper, known as a lazy quantity or lazy array, computes the corresponding output entries on-the-fly when indexed. Examples of lazy arrays in Julia include `Adjoint` arrays from the LinearAlgebra package and `SubArrays`, which represent transpositions and slices of regular Julia arrays, respectively. While eager operations are sometimes necessary – such as when solving linear systems where lazy RHS or LHS would incur prohibitive costs without advanced optimizations – lazy evaluations of cell-wise operations in Gridap, combined with Julia's JIT compilation, enable a high-level API that closely resembles the mathematical notation of weak forms of PDEs while maintaining memory efficiency and computational performance comparable to traditional codes. This approach is particularly effective because elemental operations are largely uniform across the cells of an FE mesh. By leveraging lazy arrays, the computational strategy avoids redundant allocations at the cell level, instead relying on pre-computed and reusable caches to efficiently fetch outputs when needed.

A hands-on example is provided in Listing 2. A basic familiarity with the Gridap syntax is assumed for understanding the code. In lines $7 - 12$, we define a parameter space and sample a set of parameters (with cardinality `nparams` $= 2$). Lines $14 - 25$ set up the FE mesh, FE space, and integration details using Gridap. Subsequently, we define a parameter-dependent function $\nu_p$, which is used to construct a parametric bilinear form for the stiffness matrix at line 34. The integration routine is executed in lines $34 - 35$, where the cell-wise parametric stiffness matrix `cell_mat` associated with the `Triangulation` object $\Omega$ is computed. Next, we define the connectivity structure `cell_dofs` and a parametric assembler `assem`$_p$. The assembly routine (starting from line 45) consists of three main steps:

(1) Allocate the global stiffness matrix with values initialized to zero (lines $45 - 46$).
(2) Define local cached objects for a single cell (lines $49 - 54$).

lst_param_subroutines.jl

```julia
1  using GridapROMs
2  using GridapROMs.ParamDataStructures
3  using Gridap
4  using Gridap.FESpaces, Gridap.Arrays
5
6  # Parametric space
7  pdomain = (1,5,1,5)
8  D = ParamSpace(pdomain)
9
10 # Set of parameters
11 nparams = 2
12 μ₂ = realization(D;nparams)
13
14 # Mesh
15 domain = (0,2,0,2)
16 cells = (2,2)
17 model = CartesianDiscreteModel(domain,cells)
18
19 # FE space
20 reffe = ReferenceFE(lagrangian,Float64,1)
21 V = FESpace(model,reffe)
22
23 # Integration
24 Ω = Triangulation(model)
25 dΩ = Measure(Ω,2)
26
27 # Parametric function
28 ν(μ) = x -> μ[1]*x[1]+μ[2]*x[2]
29 νₚ(μ) = parameterize(ν,μ)
30
31 # Cell-wise parametric stiffness matrix
32 v = get_fe_basis(V)
33 u = get_trial_fe_basis(V)
34 a = ∫( νₚ(μ₂)*∇(v)⋅∇(u) )dΩ
35 cell_mat = a[Ω]
36
37 # Cell-wise dof ids
38 cell_dofs = get_cell_dof_ids(V)
39
40 # Parametric assembler
41 assem = SparseMatrixAssembler(V,V)
42 assemₚ = parameterize(assem,μ)
43
44 # Allocation global parametric stiffness
45 data = ([cell_mat],[cell_dofs],[cell_dofs])
46 A = allocate_matrix(assemₚ,data)
47
48 # Allocation local caches
49 ids_cache = array_cache(cell_dofs)
50 vals_cache = array_cache(cell_mat)
51 ids1 = getindex!(ids_cache,cell_dofs,1)
52 vals1 = getindex!(vals_cache,cell_mat,1)
53 add! = FESpaces.AddEntriesMap(+)
54 add_cache = return_cache(add!,A,vals1,ids1,ids1)
55
56 # Elemental loop
57 for cell in 1:length(cell_dofs)
58     ids = getindex!(ids_cache,cell_dofs,cell)
59     vals = getindex!(vals_cache,cell_mat,cell)
60     evaluate!(add_cache,add!,A,vals,ids,ids)
61
62     # Check
63     @assert isa(vals,ParamBlock)
64 end
65
66 # Check
67 @assert isa(A,ConsecutiveParamSparseMatrixCSC)
68 @assert size(A) == (2,2)
```

FIGURE 2. Integration and assembly subroutines of a parameterized stiffness matrix. The elemental matrices returned by the integration are computed lazily, ensuring efficient memory usage and computational cost. During the assembly, memory allocation occurs only twice: first, when the global parametric stiffness matrix is initialized (line 46), and second, when caches for elemental quantities are allocated (lines $49-54$). The latter cost is minimal, as the caches are allocated once and reused across all cells in the elemental *for* loop (lines $57-64$).

(3) Perform the elemental *for* loop to assemble the global matrix from the local ones (starting from line 57).

Crucially, the last step is allocation-free, as memory consumption is confined to the previous two phases of the assembly. The significance of using lazy arrays should now be evident: they eliminate the need to allocate and compute $N_e$ elemental structures upfront, as is typical in most FE codes, where $N_e$ denotes the number of cells in the mesh. Instead, elemental values are fetched efficiently in-place (i.e., without additional allocation) one at a time, using local cached objects. Now, let us highlight the innovations introduced in GridapROMs. Due to the parameter in the bilinear form, the element type (in Julia terms, the `eltype`) of the lazy elemental stiffness matrix is no longer a standard Julia `Matrix`, as it would be in a typical Gridap application. As shown at line 63, each entry of `cell_mat` is a `ParamBlock`, representing a collection of `nparams` elemental matrices. Similarly, the output of the assembly is not a standard sparse matrix but a `ConsecutiveSparseMatrixCSC`, as illustrated at line 67. `ParamBlock` and `ConsecutiveSparseMatrixCSC` are custom types implemented in GridapROMs to represent parametric arrays at the elemental and global levels, respectively. This design offers two key advantages: first, it enables seamless reuse of Gridap's lazy and efficient implementation in a parametric context; second, it defers the *for* loop over the parameters until the global matrix entries are filled in-place, thereby avoiding unnecessary cache allocations. To clarify this point further, Fig. 3 compares the FE subroutines in GridapROMs with those implemented using a naive *for* loop over the parameters. In Fig. 4, we illustrate how the design principles of GridapROMs significantly reduce the computational cost of FE subroutines. Specifically, we compare the wall time and memory usage of GridapROMs against a naive *for* loop approach. The key distinction lies in GridapROMs' ability to reuse parametric local caches, which the naive approach neglects. As anticipated, GridapROMs demonstrates substantial computational advantages, particularly in terms of memory efficiency. Interestingly, the memory usage remains consistent across different mesh sizes, which may seem counterintuitive to those unfamiliar with lazy operations. This behavior arises because, in lazy operations, the increased size of objects primarily impacts (1) wall time and (2) the allocation of larger global and local caches. Notably, for clarity, Fig. 4 excludes the cost of global cache allocation, while the cost of local cache allocation remains negligible.

**Algorithm 1** GridapROMs subroutines.

1: Allocate $\mathtt{A}^{\mu}$
2: Compute $\mathtt{cell_mat}^{\mu}$
3: Allocate parametric local caches
4: **for** $\mathtt{cell}$ = 1:$\mathtt{\#cells}$ **do**
5: In-place fetch: $\mathtt{mat}^{\mu}$ = $\mathtt{cell_mat}^{\mu}\mathtt{[cell]}$
6: In-place fetch: $\mathtt{ids}$ = $\mathtt{cell_ids[cell]}$
7: **for** $\boldsymbol{\mu}_i \in \boldsymbol{\mu}$ **do**
8: $\mathtt{A}^{\mu_i}[\mathtt{ids},\mathtt{ids}] = \mathtt{mat}^{\mu_i}$
9: **end for**
10: **end for**

**Algorithm 2** Naive *for* loop subroutines.

1: Allocate $\mathtt{A}^{\mu}$
2: **for** $\boldsymbol{\mu}_i \in \boldsymbol{\mu}$ **do**
3: Compute $\mathtt{cell_mat}^{\mu_i}$
4: Allocate local caches
5: **for** $\mathtt{cell}$ = 1:$\mathtt{\#cells}$ **do**
6: In-place fetch: $\mathtt{mat}^{\mu_i}$ = $\mathtt{cell_mat}^{\mu_i}\mathtt{[cell]}$
7: In-place fetch: $\mathtt{ids}$ = $\mathtt{cell_ids[cell]}$
8: $\mathtt{A}^{\mu_i}[\mathtt{ids},\mathtt{ids}] = \mathtt{mat}^{\mu_i}$
9: **end for**
10: **end for**

FIGURE 3. Comparison of integration and assembly of a parametric stiffness matrix using GridapROMs (left) versus a naive outer *for* loop over the parameters (right). The GridapROMs approach is more efficient because: (1) the computation of `cell_mat` (integration) is performed simultaneously for all parameters, reusing pre-computed integration caches; (2) assembly caches are allocated once and reused across all cells; and (3) the fetching process within the elemental loop is executed only once per cell.

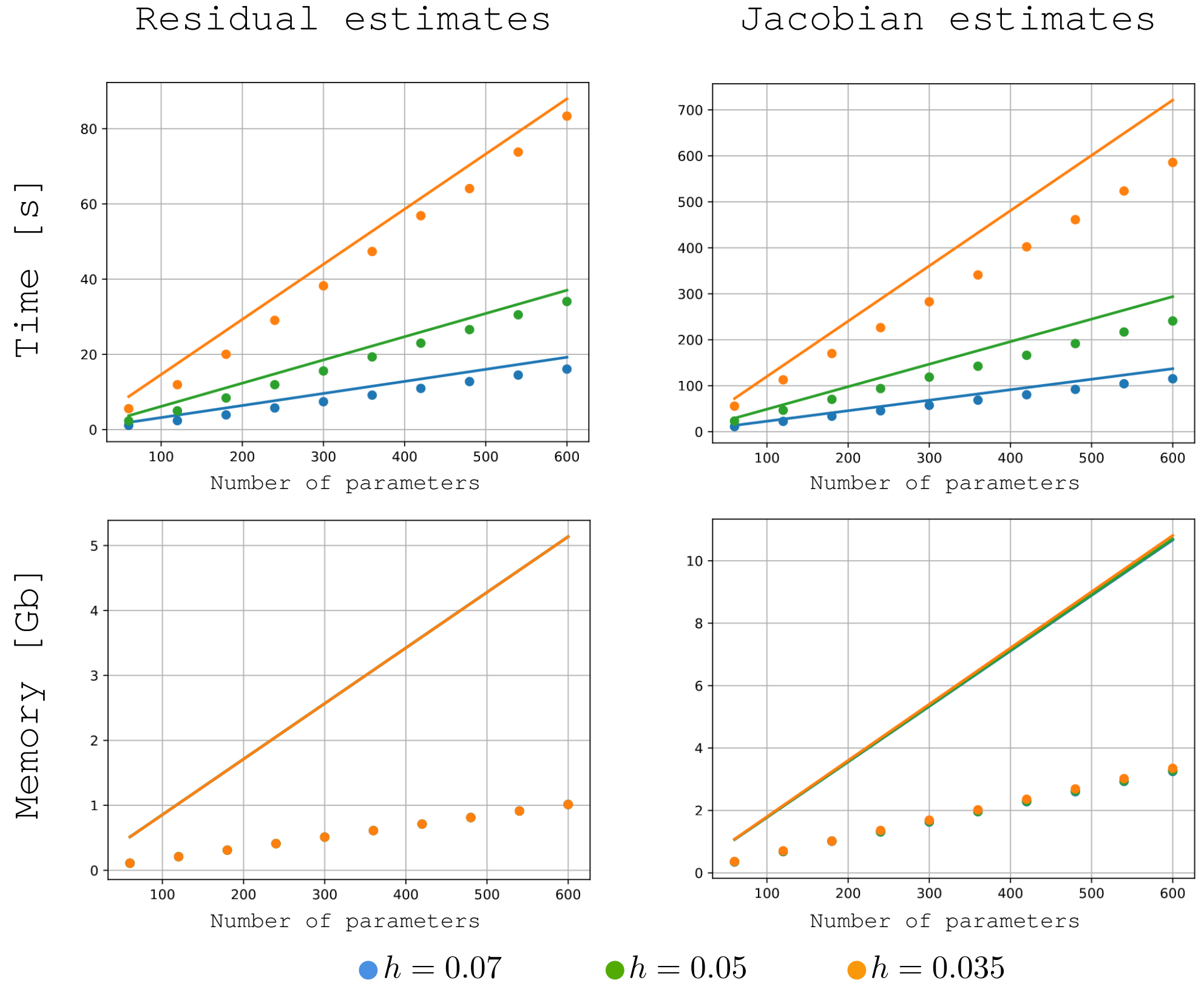

FIGURE 4. Wall time and memory usage for assembling residuals and Jacobians in GridapROMs, applied to a steady-state Navier-Stokes problem in the 3D geometry shown in Fig. 11, across varying mesh sizes. The measurements are compared against a baseline cost estimate (solid lines) as a function of the number of parameters. The baseline represents the cost of assembling a single residual/Jacobian using Gridap, scaled by the number of parameters. Note that the estimates exclude the cost of allocating global residuals and Jacobians.

**Remark 1.** *The type* `ConsecutiveSparseMatrixCSC` *represents a collection of sparse matrices with CSC ordering, where the nonzero values are stored consecutively in memory. Instead of a single vector of nonzero values, as in standard sparse matrices, it uses a matrix of nonzero values with a number of columns equal to* `nparams`*. Notably, the size of* `A` *is shown as* $(2, 2)$ *at line 68, which might seem counterintuitive since adding a parameter increases the number of entries linearly, not quadratically. However, in GridapROMs, parametric arrays are designed to maintain their original dimensionality: a parametric array of dimension* $D$ *is implemented as an array of dimension* $D$ *containing arrays of dimension* $D$*. This design choice preserves compatibility with Julia's multiple dispatch system. To accommodate this, the indexing behavior of parametric arrays with dimensions greater than* $1$ *is adjusted: accessing diagonal elements returns an*

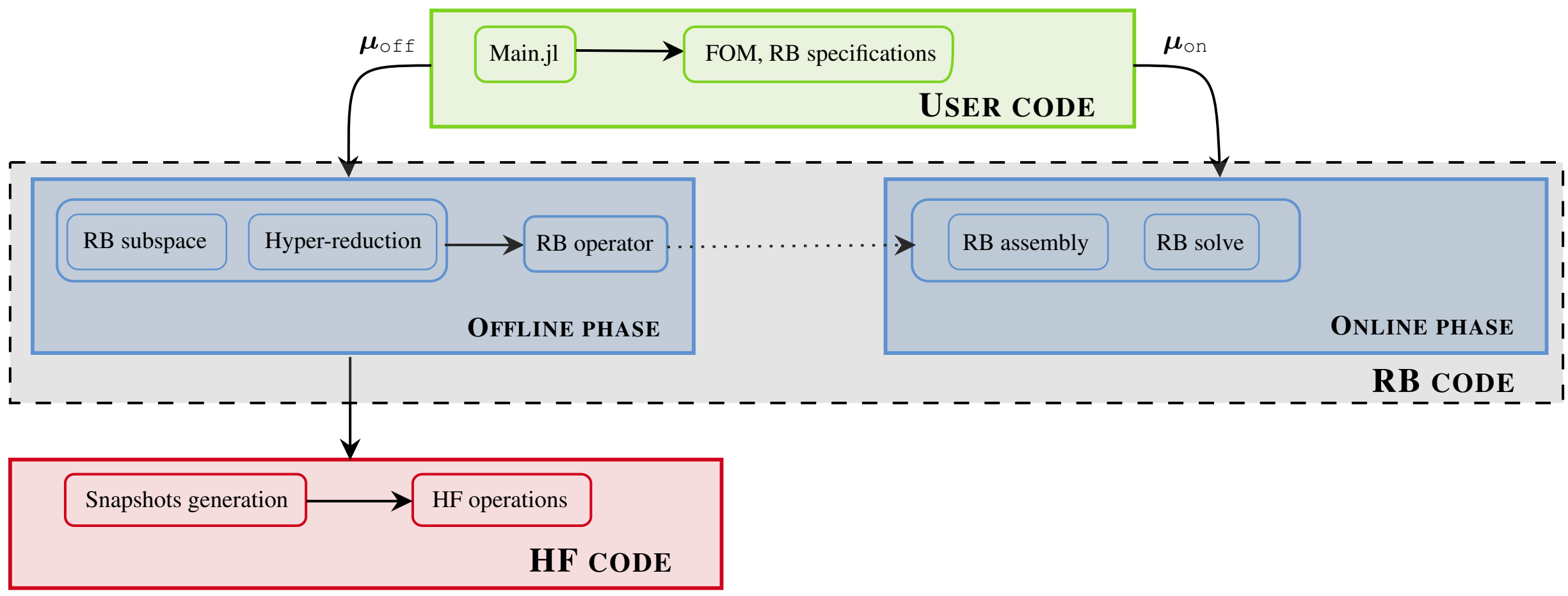

FIGURE 5. A schematic view of the implementation of GridapROMs.

| Abstract type | Purpose |
|---|---|
| `AbstractRealization` | API of realizations sampled from a parametric domain |
| `AbstractParamFunction` | Parameterized version of a Julia `Function` |
| `ParamBlock` | API of lazy, parameterized quantities defined on the FE cells |
| `ParamFEFunction` | Parameterized version of a Gridap `FEFunction` |
| `AbstractParamArray` | Parameterized version of a Julia `AbstractArray` |
| `AbstractSnapshots` | Collections of instances of `AbstractParamArray` |

TABLE 1. Main abstract types involved in the HF code in GridapROMs.

*array, while accessing off-diagonal elements returns an empty array. Thus, while a `ConsecutiveSparseMatrixCSC` is conceptually a vector of sparse matrices, it behaves like a matrix of sparse matrices for practical purposes.*

2.4. **Main abstractions.** In this subsection, we summarize the key abstractions implemented in the HF and RB codes of GridapROMs (see Fig. 5). As shown in Tb. 1, most full-order operations rely on a concise set of abstractions, which extend Gridap's lazy evaluation framework to handle parameterized problems. The first abstraction, `AbstractRealization`, represents realizations of $\mathcal{D}$. For instance, the variable defined at line 12 of Listing 2 is an `AbstractRealization`. In steady-state problems, this type acts as a wrapper for parameter sets, such as $\boldsymbol{\mu}_{\mathtt{off}}$ or $\boldsymbol{\mu}_{\mathtt{on}}$ (offline and online parameters, respectively). In transient problems, it represents tuples $\{(\boldsymbol{\mu}_i, t_j)\}_{i,j}$, where $t_j$ is the $j$th time step, enabling the computation of $(\boldsymbol{\mu}, t)$-dependent FE residuals and Jacobians in a single call[1].

The `AbstractParamFunction` abstraction provides an API for parameterized physical quantities, such as conductivity coefficients, Reynolds numbers, boundary conditions, or initial conditions. For example, $\nu_p$ in lines $27-28$ of Listing 2 is an `AbstractParamFunction` created using the `parameterize` function. During integration, an `AbstractParamFunction` is converted into a `ParamBlock`, representing parameterized quantities lazily evaluated on each cell. These evaluations yield another `ParamBlock`, containing parametric elemental vectors or matrices. These elemental quantities can then be interpolated using the FE basis, resulting in a `ParamFEFunction`.

During assembly, the list of lazy, elemental `ParamBlock` objects is converted into a global `AbstractParamArray`, which is essentially a Julia array of arrays with entries stored consecutively in memory for efficiency. The `ConsecutiveSparseMatrixCSC` type is an example of an `AbstractParamArray`, representing parameterized residuals or Jacobians. To compute solution snapshots, a *for* loop iterates over the assembled parametric arrays, solving the resulting systems of equations. This machinery at the HF code level generates the snapshots required for the offline phase. These snapshots are represented by the `AbstractSnapshots` type, which supports efficient lazy indexing, reshaping, and axis permutation.

In the RB code, instances of `AbstractSnapshots` are compressed using state-of-the-art low-rank reduction algorithms during the offline phase. The efficient indexing capabilities of `AbstractSnapshots` are crucial for the performance of these algorithms. Specifically, computationally intensive methods such as TPOD, randomized proper orthogonal decomposition (POD), and TT decompositions can be executed on `AbstractSnapshots` with efficiency comparable to that of standard Julia arrays. The outputs of these reduction algorithms are encapsulated within the `Projection` abstraction, which serves as the cornerstone of the RB code in Fig. 5. A `Projection` generally represents a mapping from a HF manifold to a RB subspace. Currently, GridapROMs supports only linear mappings, which are represented as matrices (e.g., the reduced basis $\boldsymbol{\Phi}$ introduced in (4)). However, the `Projection` abstraction is designed to accommodate nonlinear mappings, such as those based on neural networks (NNs), which are planned for future development. For the low-rank approximation of residuals and Jacobians, we use the `HyperReduction` abstraction, a specialized form of `Projection` that, in addition to the projection map, includes a reduced integration domain. This domain is essential for discrete empirical interpolation strategies to compute reduced coefficients, as detailed in Subsection 2.1. The `RBSpace` abstraction pairs a Gridap `FESpace` with the `Projection` $\boldsymbol{\Phi}$, enabling reduced solutions to be interpreted as FE functions, as described in (5). Finally, a `ReducedOperator` combines the trial and test `RBSpace` with the `HyperReduction` of residuals and Jacobians. This operator, generated during the offline phase, can be saved to a file for reuse in future simulations. Once a `ReducedOperator` is defined, the online phase can be executed, which involves assembling and solving the hyper-reduced system (11).

[1]To achieve this, modify Listing 2 as follows: (1) define a `TransientParamSpace` at line 8, including the temporal mesh; (2) define a $(\boldsymbol{\mu}, t)$-dependent function at lines $27-28$; and (3) adjust line 34 accordingly. See Subsection 2.5 for further details.

| Abstract type | Purpose |
|---|---|
| `Projection` | Projection operators from HF manifolds to RB subspaces |
| `HyperReduction` | Specialization of a `Projection` reserved for affine decompositions |
| `RBSpace` | Reduced version of a Gridap `FESpace` |
| `RBOperator` | Reduced version of a Gridap `FEOperator` |

TABLE 2. Main abstract types involved in the RB code in GridapROMs.

Lastly, we note that the low-level back-end of GridapROMs is highly customizable. Users are encouraged to implement custom types inheriting from the provided abstractions, which – thanks to the modularity of the library – allows them to retain most of the functionality already available for the existing subtypes.

2.5. **Usage example.** We provide a usage example to illustrate the expressiveness and simplicity of GridapROMs. The code in Listing 6 solves a 2D parameterized heat equation on the domain $\Omega \times [0,T] = [0,1]^2 \times [0,0.1]$, with the parameter space $\mathcal{D} = [1,5]^2$. Assuming familiarity with the Gridap API, most of the code (up to line 55) should appear as a natural extension of a standard Gridap driver for solving a heat equation. The key differences compared to a typical Gridap implementation are:

- Definition of parametric quantities, including a transient parameter space, a $(\boldsymbol{\mu}, t)$-dependent weak formulation, a trial space with a $(\boldsymbol{\mu}, t)$-dependent Dirichlet datum, and a $\boldsymbol{\mu}$-dependent initial condition.
- Introduction of triangulation lists for residuals and Jacobians, used to construct the reduced integration domains described in Subsection 2.1.
- Characterization of a parameterized `FEOperator`, incorporating the defined parametric quantities (excluding the parametric initial condition).

The RB code begins at line 58, where we define the solver `rbsolver`, which encapsulates the general RB specifications. In this example, we construct a RB from the solution snapshots using a spatio-temporal TPOD with a tolerance of $\texttt{tol} = 10^{-4}$. This method, also referred to as ST-RB [4, 5, 36], is summarized in Alg. 3 for completeness. Notably, (1) the matrix $\boldsymbol{X}$ represents a discrete inner product defined on the FE spaces of the problem, (2) the function `POD` corresponds to the standard POD [2], and (3) the mode-2 reshape at line 8 is effectively a swapping of axes rather than a conventional reshaping operation. The keyword argument `nparams`, also referred to as $N_\mu$ in Alg. 3, specifies the number of parameters

**Algorithm 3** `ST-RB`: Given a tensor of space-time snapshots $\boldsymbol{U} \in \mathbb{R}^{N \times N_t \times N_\mu}$, a prescribed accuracy `tol`, a norm matrix $\boldsymbol{X} \in \mathbb{R}^{N \times N}$, build the space-time operator $\boldsymbol{\Phi} \in \mathbb{R}^{N N_t \times n}$ that is $\boldsymbol{X}$-orthogonal in space, and $\ell^2$-orthogonal in time.

1: **function** $\texttt{ST-RB}(\boldsymbol{U}, \boldsymbol{X}, \texttt{tol})$
2: Cholesky factorization: $\boldsymbol{H}^T \boldsymbol{H} = \texttt{Cholesky}\,(\boldsymbol{X})$, ▷ $\boldsymbol{H} \in \mathbb{R}^{N \times N}$
3: Mode-1 reshape: $\boldsymbol{U}_1 = \texttt{reshape}(\boldsymbol{U}, N, N_t N_\mu)$ ▷ $\boldsymbol{U}_1 \in \mathbb{R}^{N \times N_t N_\mu}$
4: Spatial rescaling: $\widetilde{\boldsymbol{U}}_1 = \boldsymbol{H} \boldsymbol{U}_1$ ▷ $\widetilde{\boldsymbol{U}}_1 \in \mathbb{R}^{N \times N_t N_\mu}$
5: Spatial reduction: $\widetilde{\boldsymbol{\Phi}}_1 = \texttt{POD}(\widetilde{\boldsymbol{U}}_1, \texttt{tol})$ ▷ $\widetilde{\boldsymbol{\Phi}}_1 \in \mathbb{R}^{N \times n_1}$
6: Spatial inverse rescaling: $\boldsymbol{\Phi}_1 = \boldsymbol{H}^{-1} \widetilde{\boldsymbol{\Phi}}_1$ ▷ $\boldsymbol{\Phi}_1 \in \mathbb{R}^{N \times n_1}$
7: Spatial contraction: $\widehat{\boldsymbol{U}}_1 = \boldsymbol{\Phi}_1^T \boldsymbol{X} \boldsymbol{U}_1$ ▷ $\widehat{\boldsymbol{U}}_1 \in \mathbb{R}^{n_1 \times N_t N_\mu}$
8: Mode-2 reshape: $\widehat{\boldsymbol{U}}_2 = \texttt{reshape}(\widehat{\boldsymbol{U}}, N_t, n_1 N_\mu)$ ▷ $\widehat{\boldsymbol{U}}_2 \in \mathbb{R}^{N_t \times n_1 N_\mu}$
9: Temporal reduction: $\boldsymbol{\Phi}_2 = \texttt{POD}(\widehat{\boldsymbol{U}}_2, \texttt{tol})$ ▷ $\boldsymbol{\Phi}_2 \in \mathbb{R}^{N_t \times n_2}$
10: Return $\boldsymbol{\Phi} = \boldsymbol{\Phi}_1 \otimes \boldsymbol{\Phi}_2$ ▷ $\boldsymbol{\Phi} \in \mathbb{R}^{N N_t \times n}$, $n = n_1 n_2$
11: **end function**

(i.e., the number of space-time snapshots) used to construct the reduced subspace. Both `tol` and `nparams` influence

lst_heat_equation.jl

```julia
 1 using GridapROMs
 2 using Gridap
 3 using DrWatson
 4
 5 # geometry
 6 Ω = (0,1,0,1)
 7 parts = (10,10)
 8 Ωₕ = CartesianDiscreteModel(Ω,parts)
 9 τₕ = Triangulation(Ωₕ)
10
11 # temporal grid
12 θ = 0.5
13 dt = 0.01
14 t0 = 0.0
15 tf = 10*dt
16 tdomain = t0:dt:tf
17
18 # parametric quantities
19 pdomain = (1,5,1,5)
20 D  = TransientParamSpace(pdomain,tdomain)
21 u(μ,t) = x -> t*(μ[1]*x[1]^2 + μ[2]*x[2]^2)
22 uₚₜ(μ,t) = parameterize(u,μ,t)
23 f(μ,t) = x -> -Δ(u(μ,t))(x)
24 fₚₜ(μ,t) = parameterize(f,μ,t)
25
26 # numerical integration
27 order = 1
28 dΩₕ = Measure(τₕ,2order)
29
30 # weak form
31 a(μ,t,du,v,dΩₕ) = ∫(∇(v)⋅∇(du))dΩₕ
32 m(μ,t,du,v,dΩₕ) = ∫(v*du)dΩₕ
33 r(μ,t,u,v,dΩₕ) = (
34 m(μ,t,∂t(u),v,dΩₕ) + a(μ,t,u,v,dΩₕ) - ∫(fₚₜ(μ,t)*v)dΩₕ
35 )
36
37 # triangulation information
38 τₕ_a = (τₕ,)
39 τₕ_m = (τₕ,)
40 τₕ_r = (τₕ,)
41 domains = FEDomains(τₕ_r,(τₕ_a,τₕ_m))
42
43 # FE interpolation
44 reffe = ReferenceFE(lagrangian,Float64,order)
45 V = TestFESpace(Ωₕ,reffe;dirichlet_tags="boundary")
46 U = TransientTrialParamFESpace(V,uₚₜ)
47 feop = TransientParamLinearOperator((a,m),r,D,U,V,domains)
48
49 # initial condition
50 u₀(μ) = x -> 0.0
51 u₀ₚ(μ) = parameterize(u₀,μ)
52 uh₀ₚ(μ) = interpolate_everywhere(u₀ₚ(μ),U(μ,t0))
53
54 # FE solver
55 slvr = ThetaMethod(LUSolver(),dt,θ)
56
57 # RB solver
58 tol = 1e-4
59 inner_prod(u,v) = ∫(∇(v)⋅∇(u))dΩₕ
60 red_sol = TransientReduction(tol,inner_prod;nparams=20)
61 rbslvr = RBSolver(slvr,red_sol;nparams_jac=1,nparams_res=20)
62
63 dir = datadir("heat_equation")
64 create_dir(dir)
65
66 rbop = try
67 # load offline quantities
68     load_operator(dir,feop)
69 catch
70 # compute and save offline quantities
71     reduced_operator(dir,rbslvr,feop,uh₀ₚ)
72 end
73
74 # online phase
75 μₒₙ = realization(feop;nparams=10,sampling=:uniform)
76 x̂,rbstats = solve(rbslvr,rbop,μₒₙ,uh₀ₚ)
77
78 # post process
79 x,stats = solution_snapshots(slvr,feop,μₒₙ,uh₀ₚ)
80 perf = eval_performance(rbslvr,feop,rbop,x,x̂,stats,rbstats)
```

FIGURE 6. Solving a parameterized heat equation with GridapROMs.

```
julia> perf
"----------------- RBPerformance ------------------
> error: 7.814050983154542e-5
> speedup in time: 55.64980423839414
> speedup in memory: 25.51763919028245
---------------------------------------------------"
```

FIGURE 7. Parameterized heat equation results.

the quality of the ROM approximation and, consequently, the accuracy of the method. Therefore, these hyperparameters should be carefully selected. Typically, `tol` is chosen within the range $[10^{-5}, 10^{-1}]$, while `nparams` is determined based on `tol`. Ideally, we aim to minimize the cost of snapshots computation by selecting a small `nparams`, such as $\sim 1$. However, the snapshots must adequately represent the manifolds being approximated (e.g., the solution manifold and those of residuals/Jacobians during hyper-reduction), necessitating an appropriate number of parameters sampled from $\mathcal{D}$. Since both the accuracy and computational cost increase with `nparams`, it is crucial to balance these factors to achieve an optimal cost-benefit ratio. For instance, a small `nparams` suffices for large tolerances $(10^{-2} - 10^{0})$, whereas `nparams` should increase for smaller `tol` values to improve accuracy. Additionally, `nparams` should be increased for more complex manifolds, such as those arising in nonlinear applications. In simpler cases, like the usage example, good accuracy can be achieved with lower `nparams` values.

To define the hyper-reduction strategy for residuals and Jacobians, the keyword arguments `nparams_res` and `nparams_jac` are passed when defining `rbsolver`. Note that no hyper-reduction is needed for $\boldsymbol{\mu}$-independent

Jacobians; in such cases, setting `nparams_jac` $= 1$ suffices. Next, we attempt to load the `ReducedOperator` from a file; if unavailable (e.g., during the first run), the offline phase is executed. Once a `ReducedOperator` is obtained, the online phase can be run for any set of online realizations $\boldsymbol{\mu}_{\text{on}}$ (disjoint from the offline parameters). This set can be generated using the `realization` function from the parameter space or, as in this example, from the FE operator. The keyword `uniform` ensures that parameters are uniformly distributed over $\mathcal{D}$; otherwise, the default sampling uses a Halton sequence [37], which provides better coverage of the sampling space compared to uniform distribution. GridapROMs also supports other sampling strategies, such as normal distribution, Latin Hypercube sampling [38], and tensorial uniform sampling [2].

Finally, the algorithm's performance relative to HF simulations can be evaluated. This involves collecting HF solutions for $\boldsymbol{\mu}_{\text{on}}$. The final call to `eval_performance` returns:

- The relative error between the HF and RB solutions is computed by averaging over all values in $\boldsymbol{\mu}_{\text{on}}$. In steady problems, this error is measured using the norm specified by `inner_prod`, which in this case corresponds to the $H_0^1$ norm. In contrast, for time-dependent applications, the spatial error is evaluated at each time step and then integrated over the entire time interval, resulting in the following error measure:

$$\frac{1}{|\boldsymbol{\mu}_{\text{on}}|} \sum_{\mu \in \boldsymbol{\mu}_{\text{on}}} \frac{\sum_{n=1}^{N_t} \|\boldsymbol{u}_{(n)}^{\mu} - (\boldsymbol{\Phi}\widehat{\boldsymbol{u}}^{\mu})_{(n)}\|_{\texttt{inner_prod}}}{\sum_{n=1}^{N_t} \|\boldsymbol{u}_{(n)}^{\mu}\|_{\texttt{inner_prod}}}.$$

- The computational speedup achieved by the RB online code compared to the HF simulations, measured in terms of wall time and memory usage. The speedup is calculated as the ratio of the HF cost metrics in `stats` to the corresponding RB metrics in `rbstats`.

The results of `eval_performance` are presented in Listing 7.

## 3. Advanced utilities

In this section, we highlight some of the more advanced utilities provided by the library. Up to this point, we have shown how a ST-RB subspace can be constructed to approximate a simple parabolic equation. We now turn to tools available in GridapROMs that extend these ideas to more challenging settings – such as saddle-point problems [39] and equations posed on moving domains. We also discuss how distributed-in-memory parallelism can be leveraged to accelerate snapshots generation and improve offline computational efficiency.

### 3.1. **Supremizer enrichment.**

When constructing a ROM for a saddle-point problem, it is often necessary to explicitly enforce inf-sup stability through a suitable stabilization strategy. A widely used approach is the supremizer enrichment of the velocity RB subspace [4, 40–43]. The procedure is outlined in Alg. 4 for the steady case; for details on its extension to unsteady saddle-point problems, we refer the reader to [4].

**Algorithm 4** ENRICHMENT: Given the velocity basis $\boldsymbol{\Phi}^u \in \mathbb{R}^{N^u \times n^u}$, the pressure basis $\boldsymbol{\Phi}^p \in \mathbb{R}^{N^p \times n^p}$, the norm matrix $\boldsymbol{X} \in \mathbb{R}^{N^u \times N^u}$ representing the inner product on the FE velocity space, and the velocity-pressure coupling matrix $\boldsymbol{B} \in \mathbb{R}^{N^p \times N^u}$, return the $N^u \times (n^u + n^p)$-dimensional enriched velocity basis.

1: **function** ENRICHMENT($\boldsymbol{\Phi}^u, \boldsymbol{\Phi}^p, \boldsymbol{X}, \boldsymbol{B}$)
2: Cholesky factorization: $\boldsymbol{H}^T\boldsymbol{H} = \texttt{Cholesky}\left(\boldsymbol{X}^u\right)$ ▷ $\boldsymbol{H} \in \mathbb{R}^{N^u \times N^u}$
3: Compute supremizers: $\boldsymbol{S} = \boldsymbol{H}^{-1}\boldsymbol{B}^T\boldsymbol{\Phi}^p$ ▷ $\boldsymbol{S} \in \mathbb{R}^{N^u \times n^p}$
4: Supremizer enrichment: $\boldsymbol{\Phi}^u \leftarrow \texttt{Gram-Schmidt}\left([\boldsymbol{\Phi}^u, \boldsymbol{S}], \boldsymbol{X}\right)$ ▷ $\boldsymbol{\Phi}^u \in \mathbb{R}^{N^u \times (n^u + n^p)}$
5: **return** $\boldsymbol{\Phi}^u$
6: **end function**

Applying this enrichment strategy yields an inf-sup stable ROM [41]. In GridapROMs, the method is available for both steady and transient saddle-point problems. Its use is straightforward: at the user level, enrichment is activated by adding a single line of code, as illustrated in Listing 8. The only additional input required is the bilinear form defining the velocity-pressure coupling, which is passed to the `Reduction` constructor. From that point, GridapROMs automatically performs the supremizer enrichment once both velocity and pressure snapshots have been compressed. A complete example script demonstrating the solution of a saddle-point problem with GridapROMs is provided in Appendix A.

### 3.2. **Tensor-train reduction.**

In this section, we describe how the TT-RB procedure introduced in [27] can be implemented. Without delving into unnecessary technical detail, we recall that this strategy significantly reduces the computational cost of constructing RBs that are orthogonal with respect to an energy norm (no such benefit arises when using Euclidean-orthogonal bases). It also provides a cleaner framework for compressing high-dimensional snapshots, such as those arising in transient simulations. While proposed only for problems defined on Cartesian grids, the subsequent work [44] demonstrates how this approach can be generalized to arbitrary domains, including cases with shape parameterizations.

lst_supremizer_enrichment.jl

```
tol = 1e-4
nparams = 50
energy((du,dp),(v,q)) = ∫(∇(du)⊙∇(v))Ωₕ + ∫(dp*q)Ωₕ
coupling((du,dp),(v,q)) = ∫(dp*(∇⋅(v)))Ωₕ
red_sol = Reduction(coupling,tol,energy;nparams)
# For transient applications:
# red_sol = TransientReduction(coupling,tol,energy;nparams)
```

FIGURE 8. Defining a reduction strategy supporting the supremizer enrichment of the velocity basis with GridapROMs.

As with supremizer enrichment, only minor modifications are required to activate this reduction type. Consider the example in Subsection 2.5. To solve the heat equation in Listing 6 using TT-RB instead of ST-RB, we only need to make the following adjustments:

- **Discrete model**: Define $\Omega_h$ as a `TProductDiscreteModel` instead of a `CartesianDiscreteModel`, using the same input arguments as in line 8 of Listing 6. Unlike a standard Gridap discrete model, this type also stores information about the one-dimensional discrete models whose tensor product defines $\Omega_h$. This only applies when the underlying domain is Cartesian. Thanks to multiple dispatch, corresponding test and trial FE spaces are created that likewise store the associated one-dimensional FE spaces. In practice, this enables the assembly of FE matrices and vectors by exploiting tensor-product decompositions, so that snapshots are represented as generalized multi-dimensional tensors rather than as conventional concatenations of matrices or vectors.
- **Tolerance vector**: Instead of a scalar tolerance, define `tol` as a vector of three values. The rationale is straightforward: in this application, the snapshots are viewed as a four-dimensional array, with the first two axes corresponding to spatial entries, the third to temporal values, and the fourth to parameters. A tensor-based rank-reduction algorithm – currently, only TT-SVD is implemented – is then applied to compress the first three axes (spatial and temporal). More generally, if snapshots are stored as a D-dimensional tensor, `tol` should be a vector of length D-1, and TT-SVD will compress the first D-1 dimensions of the tensor.

For problems posed on non-Cartesian domains, the TT-RB approach can still be employed by combining it with an unfitted FE discretization [45–49] – which relies on a background Cartesian mesh. A detailed discussion of this methodology can be found in [44]. For completeness, Appendix B provides the code used to solve a Poisson equation on a parameterized domain using TT-RB.

3.3. **Local subspace strategy.** Local reduced subspaces have been introduced to improve reducibility in applications where standard global subspaces fail to provide meaningful dimensionality reduction [22–24]. In GridapROMs, we adopt the local subspace strategy presented in [44], which builds on the earlier work of [26]. The methodology is outlined in Algs. 5-6. For additional details, we refer to [26, 44].

**Algorithm 5** `OFFLINE`: Given a set of offline parameters $\boldsymbol{\mu}_1, \ldots, \boldsymbol{\mu}_{N_\mu} \in \mathcal{D}$, the number of clusters for the reduced subspace $N_c$, and the number of clusters for the hyper-reduction $N_c^h$, compute all offline quantities.

1: **function** `OFFLINE`($\boldsymbol{\mu}_1, \ldots, \boldsymbol{\mu}_{N_\mu}, N_c, N_c^h$)
2: Run clustering of parameters (reduced subspace): $\boldsymbol{\alpha}_1, \ldots, \boldsymbol{\alpha}_{N_c} = \texttt{Kmeans}([\boldsymbol{\mu}_1, \ldots, \boldsymbol{\mu}_{N_\mu}], N_c)$,
3: Run clustering of parameters (hyper-reduction): $\boldsymbol{\beta}_1, \ldots, \boldsymbol{\beta}_{N_c^h} = \texttt{Kmeans}([\boldsymbol{\mu}_1, \ldots, \boldsymbol{\mu}_{N_\mu}], N_c^h)$
4: Compute solution, LHS and RHS snapshots for every $\boldsymbol{\mu}_1, \ldots, \boldsymbol{\mu}_{N_\mu}$
5: **for** $j = 1, \ldots, N_c$ **do**
6: Compute the $j$th RB subspace from the solution snapshots associated with $\boldsymbol{\alpha}_j$
7: **for** $k = 1, \ldots, N_c^h$ **do**
8: Compute the $k$th hyper-reductions from the residual and Jacobian snapshots associated with $\boldsymbol{\beta}_k$
9: Project the $k$th hyper-reductions onto the $j$th RB subspace
10: **end for**
11: **end for**
12: **return** the $N_c$ subspaces and $N_c \cdot N_c^h$ Galerkin projections
13: **end function**

From the user's perspective, enabling localization requires only minimal changes compared to the default construction of global subspaces and hyper-reduction structures. As shown in Listing 9, the localized strategy is activated simply by specifying an appropriate reduction method for both solutions and residuals/Jacobians. A complete usage example demonstrating this localization procedure is provided in Appendix B.

3.4. **Distributed-in-memory snapshots.** Lastly, we illustrate how GridapROMs can exploit parallel computing to accelerate the generation of HF snapshots. Parallelism is achieved through the following Julia packages:

**Algorithm 6** `ONLINE`: Given an online parameter $\boldsymbol{\mu} \in \mathcal{D}$ and the quantities computed offline, compute the reduced solution $\boldsymbol{u}_n(\boldsymbol{\mu}) \in \mathbb{R}^n$.

1: **function** `ONLINE`($\boldsymbol{\mu}$)
2: Find closest subspace cluster: $j = \arg \min_{i=1,\dots,N_c} \|\boldsymbol{\mu} - \boldsymbol{\alpha}_i\|_2$
3: Find closest hyper-reduction cluster: $k = \arg \min_{i=1,\dots,N_c^h} \|\boldsymbol{\mu} - \boldsymbol{\beta}_i\|_2$
4: Fetch the offline quantities associated with $(j, k)$
5: Run the hyper-reduction's online phase (see (10))
6: Assemble and solve the hyper-reduced system (see (11)-(12))
7: **return** $\boldsymbol{u}_n^{\mu}$
8: **end function**

lst_localization.jl

```
1  tol = 1e-4
2  inner_prod(du,v) = ∫(∇(du)⋅∇(v))Ωₕ
3
4  kwargs = (
5    nparams = 50,
6    compression = :local,
7    ncentroids = 5
8  )
9
10 slvr = LUSolver()
11 red_sol = Reduction(tol,inner_prod;kwargs...)
12 red_res = HyperReduction(tol;kwargs...)
13 red_jac = HyperReduction(tol;kwargs...)
14 rbslvr = RBSolver(slvr,red_sol,red_res,red_jac)
15
16 # More concisely:
17 nparams_res = nparams_jac = nparams
18 rbslvr = RBSolver(slvr,red_sol;nparams_res,nparams_jac)
19
20 # For transient applications:
21 # dt,θ = ...
22 # odeslvr = ThetaMethod(slvr,dt,θ)
23 # red_sol = TransientReduction(tol,inner_prod;kwargs...)
24 # rbslvr = RBSolver(odeslvr,red_sol;nparams_res,nparams_jac)
```

FIGURE 9. Defining a reduction strategy supporting the localization of subspaces and hyper-reduction quantities with GridapROMs.

- GridapDistributed [50]: extends Gridap to support the parallel FE solution of PDEs. It builds on PartitionedArrays [51], which provides distributed linear algebra operations.
- GridapPETSc: offers an interface to Portable, Extensible Toolkit for Scientific Computation (PETSc) [52], enabling the use of a broad suite of highly scalable solvers within Julia.
- GridapSolvers [53]: complements GridapPETSc with physics-informed solvers, including Krylov-based iterative methods and geometric multigrid algorithms.

These packages form a robust and efficient ecosystem for HPC applications in Julia. They have been successfully applied to challenging large-scale problems, achieving excellent strong and weak scaling results – for example, in the steady Stokes system and the nonlinear rotating shallow-water equations on the sphere. GridapROMs extends this parallelism to efficiently handle parameterized PDEs, following the same design principles as in the serial case. At present, distributed capabilities are available for the generation of solution, residual, and Jacobian snapshots; future work will extend this machinery to support all offline computations.

In Listing 10, we provide a usage example where 10 snapshots are computed by solving a parametric Poisson equation, distributing the computations across 4 processors. Since this script cannot be run interactively, it must be executed from the terminal using:

```
mpiexecjl --project=.  -n 4 julia lst_distributed_snapshots.jl
```

The script differs only slightly from its serial counterpart. The most notable change is that the execution of the `main` function is wrapped inside a Julia `do`-block (lines 69-72). Two arguments are passed to `main`: (i) the process grid layout, represented as a tuple of process ranks, and (ii) a `distribute` function that creates a distributed array holding these ranks. Each of the four processors receives exactly one rank, i.e., one partition.

Another key difference is the definition of a second `do`-block, which wraps most of `main` (lines 17-63). Here, `GridapPETSc.with` is called with two arguments: (i) a zero-argument function whose body (essentially, lines 18-62) consists of the instructions to be executed on each process, and (ii) a string of PETSc solver options. For a comprehensive explanation of these options, we refer to the official PETSc documentation [52]. In this example, we use the following configuration:

- `-ksp_type cg`: use the conjugated gradient method as the Krylov subspace solver.
- `-pc_type gamg`: use a geometric algebraic multigrid preconditioner.
- `-ksp_monitor`: print the solver progress at each iteration.
- `-ksp_rtol 1e-8`: terminate the iterations once the relative error falls below a tolerance of $10^{-8}$.

Finally, we note that these options need not be explicitly passed to the `PETScLinearSolver` constructor.

lst_distributed_snapshots.jl

```julia
using Gridap
using GridapDistributed
using GridapROMs
using GridapPETSc
using PartitionedArrays
using Test

function main(distribute,parts)
  # Distributed array with the identifiers
  # of the parallel processes
  ranks = distribute(LinearIndices((prod(parts),)))

  # PETSc linear solver options
  options = "-ksp_type cg -pc_type gamg -ksp_monitor
             -ksp_rtol 1e-8"

  GridapPETSc.with(args=split(options)) do
    # geometry
    domain = (0,1,0,1)
    cells = (40,40)
    Ωₕ = CartesianDiscreteModel(ranks,parts,domain,cells)
    τₕ = Triangulation(Ωₕ)

    # parametric quantities
    D = ParamSpace((0,1))
    u(μ) = x -> μ[1]*x[1]+x[2]
    f(μ) = x -> -Δ(u(μ))(x)
    fₚ(μ) = parameterize(f,μ)
    uₚ(μ) = parameterize(u,μ)

    # numerical integration
    order = 2
    degree = 2*order
    dΩₕ = Measure(τₕ,degree)

    # weak form
    a(μ,u,v) = ∫( ∇(v)⋅∇(u) )dΩₕ
    l(μ,u,v) = a(μ,u,v) - ∫( v*fₚ(μ) )dΩₕ

    # FE interpolation
    reffe = ReferenceFE(lagrangian,Float64,order)
    V = TestFESpace(Ωₕ,reffe,dirichlet_tags="boundary")
    U = ParamTrialFESpace(V,uₚ)
    feop = LinearParamOperator(l,a,D,U,V)

    # FE solver
    slvr = PETScLinearSolver()

    # FE snapshots
    μ = realization(D;nparams=10)
    sol, = solve(slvr,feop,μ)

    # Post process
    for (i,μᵢ) in enumerate(μ)
      solᵢ = param_getindex(sol,i)

      Uᵢ = TrialFESpace(V,u(μᵢ))
      uₕᵢ = FEFunction(Uᵢ,solᵢ)

      eₕᵢ = u(μᵢ) - uₕᵢ
      @test sum( ∫(eₕᵢ*eₕᵢ)dΩₕ ) < 1e-8
    end
  end
end

# process grid layout
parts = (2,2)

# launch main in parallel
with_mpi() do distribute
  main(distribute,parts)
end
```

FIGURE 10. Generating 10 snapshots with values distributed across 4 processors in GridapROMs.

## 4. APPLICATION

In this section, we present the numerical results obtained using GridapROMs to solve a fluid dynamics problem modeled by the unsteady Navier-Stokes equations (16) in a 3D geometry, as illustrated in Fig. 11.

$$\begin{cases} \frac{d\underline{u}^\mu}{dt} + \nabla \cdot (\nu^\mu \nabla \underline{u}^\mu) + (\underline{u}^\mu \cdot \nabla)\underline{u}^\mu - \nabla p^\mu = \underline{0} & (\underline{x}, t) \in \Omega \times (0, T] \\ \nabla \cdot \underline{u}^\mu = 0 & (\underline{x}, t) \in \Omega \times (0, T]; \\ \underline{u}^\mu = \underline{g}^\mu & (\underline{x}, t) \in \Gamma_D \times (0, T]; \\ \nu^\mu \nabla \underline{u}^\mu \cdot \underline{n} - p\underline{n} = \underline{0} & (\underline{x}, t) \in \Gamma_N \times (0, T]; \\ \underline{u}^\mu = \underline{0} & (\underline{x}, t) \in \Omega \times \{0\}. \end{cases} \tag{16}$$

The domain $\Omega$ is a rectangular prism with dimensions $(L, W, H) = (1, 0.5, 0.1)$, featuring two cylindrical holes of radius $R = 0.1$ and height $H$. The problem involves a viscosity and inflow that vary with both time and parameters, modeled by a Dirichlet condition on the inlet boundary. A homogeneous Neumann condition is imposed on the outlet, while the remaining walls are subject to a homogeneous Dirichlet condition (with flow constrained only in the normal direction on the top and bottom walls). The initial condition is also homogeneous. The parametric data are defined as follows:

$$\nu^\mu(\underline{x}, t) = \frac{\mu_1}{100}; \qquad \underline{g}^\mu(\underline{x}, t) = -x_2(W - x_2)\left(1 - \cos(\pi t/T) + \frac{\mu_3}{100}\sin(\mu_2 \pi t/T)\right)\underline{n}_1,$$

where $\underline{n}_1 = (1, 0, 0)^T$. The temporal domain is $[0, 0.15]$, discretized into 60 uniform time steps, and the parameter space is defined as $\mathcal{D} = [1, 10]^3$. For the spatial discretization, we employ the inf-sup stable pair of FE spaces $(\mathcal{V}_h, \mathcal{Q}_h) = (P_2, P_1)$ for the velocity and pressure, respectively. The spatial DOFs count is 15943 for the velocity and 1211 for the pressure, resulting in a total of $N = 1029240$ space-time DOFs. For temporal discretization, we use the BE time-marching scheme. The solution snapshots are generated on the `GADI`[2] supercomputer, where the FE code is executed on 10 processors, each handling 6 parameters across all 60 time steps. The snapshots are then concatenated into a single dataset. Out of these, 55 snapshots are used to construct a $(H^1)^3$-orthogonal RB for the velocity and an $L^2$-orthogonal RB for the pressure using a Sparse Random Gaussian technique [21], while the remaining 5 snapshots form the test set. Additionally, we apply an inf-sup stabilization procedure via supremizer enrichment of the velocity basis [4, 40–43], a standard approach for RB approximations of saddle-point problems [39], such as the Navier-Stokes equations. For hyper-reduction, we use a space-time MDEIM technique [5], with `nparams_res` = 55 for the residual and `nparams_jac` = 15 for the Jacobian. Results for tolerances `tol` $\in \{10^{-i}\}_{i=3}^{5}$, commonly used to evaluate ROMs, are presented in Fig. 12 and Tb. 3. In particular, we are interested in evaluating the error measure

$$\mathbf{E} = \begin{bmatrix} \sum\limits_{\boldsymbol{\mu}\in\boldsymbol{\mu}_{\text{on}}} \left( \int_0^T \left( \|\boldsymbol{u}_h^\mu(t) - \boldsymbol{u}_n^\mu(t)\|_{(H^1(\Omega))^3} / \|\boldsymbol{u}_h^\mu(t)\|_{(H^1(\Omega))^3} \right) dt \right) / N_{\text{on}} \\ \sum\limits_{\boldsymbol{\mu}\in\boldsymbol{\mu}_{\text{on}}} \left( \int_0^T \left( \|\boldsymbol{p}_h^\mu(t) - \boldsymbol{p}_n^\mu(t)\|_{L^2(\Omega)} / \|\boldsymbol{p}_h^\mu(t)\|_{L^2(\Omega)} \right) dt \right) / N_{\text{on}} \end{bmatrix},$$

where $N_{\text{on}} = 5$ represents the number of online parameters. Additionally, we evaluate the computational speedup in time, **SU-T**, defined as the ratio of the average wall time for a FOM simulation to that of a ROM simulation. Similarly, the memory speedup, **SU-M**, is defined as the ratio of the average memory allocations in a FOM simulation to those in a ROM simulation.

[2] https://nci.org.au/

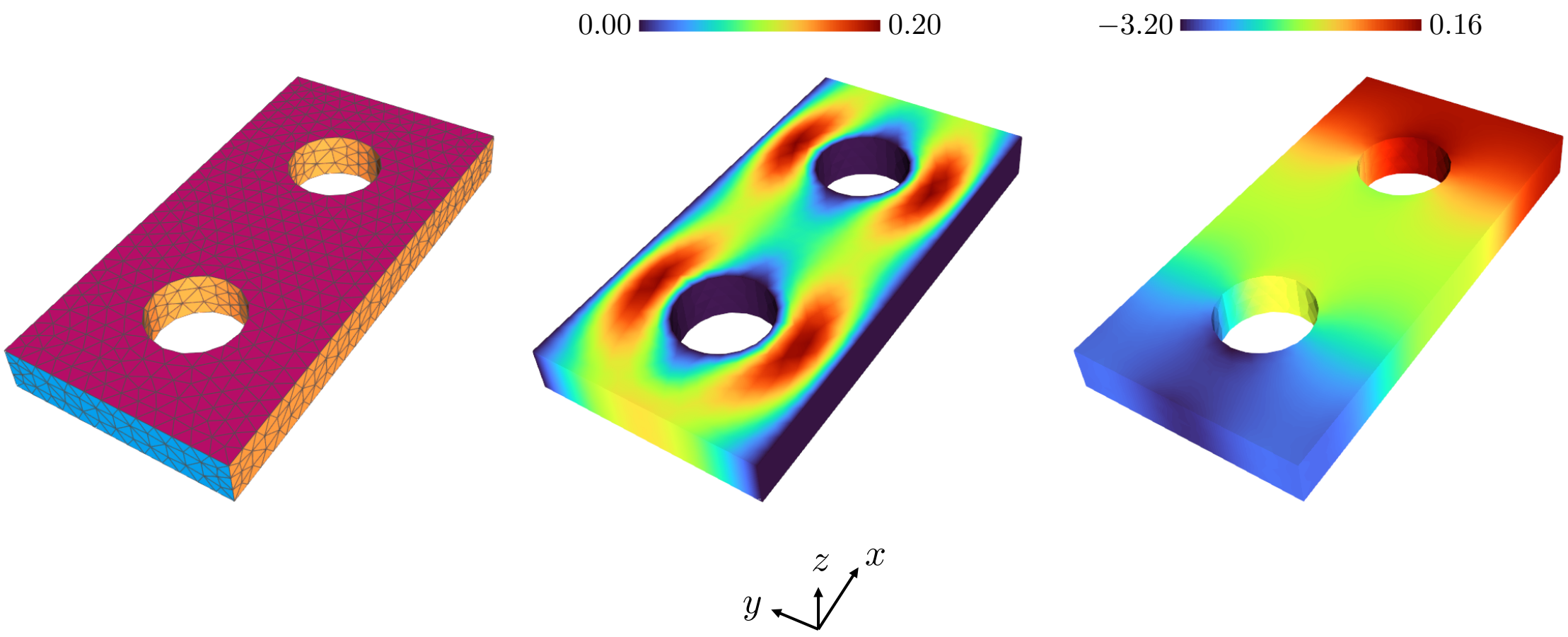

FIGURE 11. Geometry used for the numerical test case (left), and the FE solution (velocity magnitude in the middle, and pressure on the right) for the parameter $\boldsymbol{\mu} = (7.35, 2.00, 4.48)$ at $t = T$. On the side walls and cylinders (orange, left figure), a no-slip Dirichlet condition is imposed; on the inlet (blue), a non-homogeneous Dirichlet condition is applied; on the top and bottom walls (magenta), a no-penetration Dirichlet condition is enforced; and on the outlet (not shown, opposite the inlet), a homogeneous Neumann condition is imposed. The union of the side walls, inlet, and top and bottom facets is denoted as $\Gamma_D$, while the outlet is denoted as $\Gamma_N$.

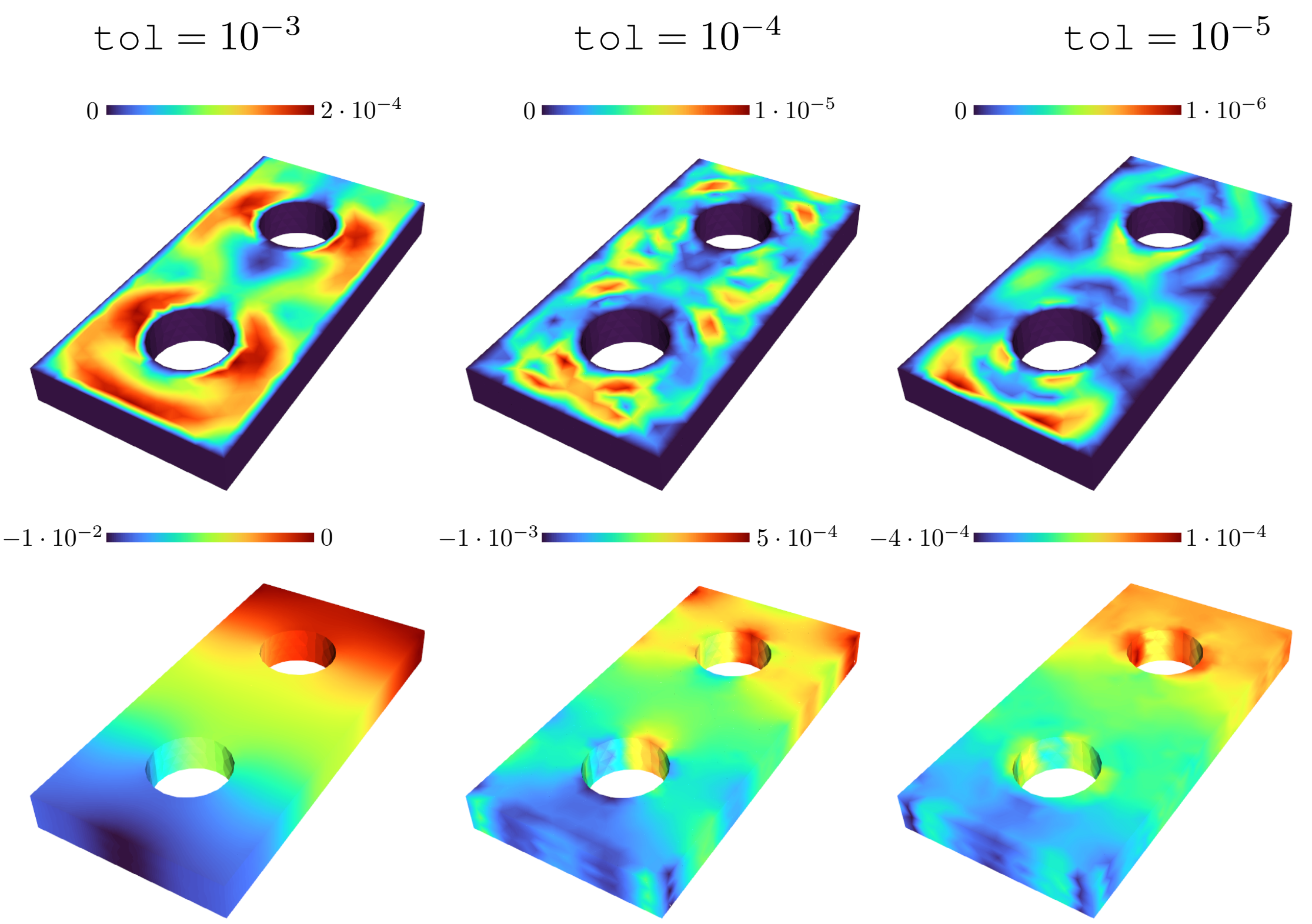

FIGURE 12. Point-wise error between the FE solution for the parameter $\boldsymbol{\mu} = (7.35, 2.00, 4.48)$ at $t = T$ and the corresponding solutions computed with GridapROMs for various tolerances. The first row shows the velocity magnitude errors for tolerances $\{10^{-i}\}_{i=3}^{5}$, while the second row presents the pressure field errors.

| $\texttt{tol} = 10^{-3}$ | | | | $\texttt{tol} = 10^{-4}$ | | | | $\texttt{tol} = 10^{-5}$ | | | |
|---|---|---|---|---|---|---|---|---|---|---|---|
| $\boldsymbol{n}$ | **E** | **SU-T** | **SU-M** | $\boldsymbol{n}$ | **E** | **SU-T** | **SU-M** | $\boldsymbol{n}$ | **E** | **SU-T** | **SU-M** |
| $\begin{bmatrix} 96 \\ 24 \end{bmatrix}$ | $\begin{bmatrix} 3 \cdot 10^{-3} \\ 2 \cdot 10^{-2} \end{bmatrix}$ | $3 \cdot 10^5$ | $7 \cdot 10^2$ | $\begin{bmatrix} 182 \\ 56 \end{bmatrix}$ | $\begin{bmatrix} 1 \cdot 10^{-3} \\ 1 \cdot 10^{-2} \end{bmatrix}$ | $3 \cdot 10^5$ | $6 \cdot 10^2$ | $\begin{bmatrix} 336 \\ 96 \end{bmatrix}$ | $\begin{bmatrix} 3 \cdot 10^{-4} \\ 1 \cdot 10^{-3} \end{bmatrix}$ | $2 \cdot 10^5$ | $4 \cdot 10^2$ |

TABLE 3. From left to right: dimensions of the reduced subspaces for velocity and pressure, average relative space-time errors for velocity and pressure, average computational speedup in execution time, and average computational speedup in memory usage, for various tolerance values.

Our findings in Tb. 3 highlight the ROM's capability to deliver highly accurate solutions at a fraction of the computational cost of FE simulations. The significant speedup, particularly in execution time, stems from the efficient implementation of the HF code with a minimal memory footprint. In terms of accuracy, the reduced subspace dimension increases as the tolerance decreases, leading to progressively more accurate ROM solutions. This trend is clearly illustrated in Fig. 12, where the point-wise errors diminish consistently with smaller `tol` values.

## 5. CONCLUSIONS AND FUTURE WORK

In this work, we introduce GridapROMs, a Julia-based library designed for solving parameterized PDEs using ROMs. By combining a user-friendly high-level API, the performance benefits of Julia's JIT compiler, and extensive use of lazy evaluations, the library achieves both flexibility and efficiency. GridapROMs supports a broad spectrum of applications, including steady, transient, single-field, multi-field, linear, and nonlinear problems. The library also provides advanced utilities, such as stabilization techniques for ROMs applied to saddle-point problems, TT-based ROMs, localization strategies, and distributed-in-memory machinery for snapshots generation, all of which are documented with illustrative code snippets. We further demonstrate its capabilities by solving a fluid dynamics problem modeled by the transient Navier-Stokes equations in a 3D geometry, showcasing significant computational cost reductions compared to HF simulations while maintaining high accuracy. Furthermore, benchmarks involving shape parameterizations are also implemented and discussed in the appendices. Additional numerical examples – including linear elasticity, electromagnetic wave scattering, and advection-diffusion equations – are provided and documented on the GridapROMs website.

We foresee several directions for advancing GridapROMs. First, we aim to extend the distributed-in-memory framework beyond snapshots generation to encompass all offline computations, including snapshots compression, hyper-reduction, and Galerkin projections. Second, we plan to benchmark the performance of GridapROMs against that of other well-established ROM libraries, both in serial and parallel settings. Third, we intend to enrich the set of reduction methods by incorporating subspaces constructed through the greedy selection of snapshots. Fourth, we aim to integrate nonlinear, deep learning-based models into the framework. Nonlinear ROMs leveraging autoencoder architectures have demonstrated the ability to construct lower-dimensional latent spaces that capture nonlinear dynamics more effectively than linear compression methods such as POD in several applications [31, 54]. Fifth, we plan to explore spatially localized strategies, made possible by distributing the FE mesh across multiple processors. Recent studies show that shallow autoencoder architectures combined with domain-decomposition pipelines can achieve both high accuracy and strong scalability [55]. In particular, it would be worthwhile to investigate how GridapROMs' existing parameter-space localization strategies could be integrated with such spatial localization approaches.

## ACKNOWLEDGMENTS

This research was partially funded by the Australian Government through the Australian Research Council (project numbers DP210103092 and DP220103160). Computational resources were provided by the Australian Government through NCI under the NCMAS and by Monash through the by Monash-NCI scheme for HPC services. N. Mueller acknowledges the Monash Graduate Scholarship from Monash University.

## REFERENCES

[1] T. Lassila, A. Manzoni, A. Quarteroni, and G. Rozza. "Model order reduction in fluid dynamics: challenges and perspectives". In: *MATHICSE-CMCS Modelling and Scientific Computing* (2013).

[2] A. Quarteroni, A. Manzoni, and F. Negri. *Reduced basis methods for partial differential equations: an introduction*. Vol. 92. Springer, 2015.

[3] F. Negri, A. Manzoni, and D. Amsallem. "Efficient model reduction of parametrized systems by matrix discrete empirical interpolation". In: *Journal of Computational Physics* 303 (2015), pp. 431–454. DOI: https://doi.org/10.1016/j.jcp.2015.09.046.

[4] R. Tenderini, N. Mueller, and S. Deparis. "Space-Time Reduced Basis Methods for Parametrized Unsteady Stokes Equations". In: *SIAM Journal on Scientific Computing* 46.1 (2024), B1–B32. DOI: 10.1137/22M1509114. eprint: https://doi.org/10.1137/22M1509114.

[5] N. Mueller and S. Badia. "Model order reduction with novel discrete empirical interpolation methods in space-time". In: *Journal of Computational and Applied Mathematics* 444 (2024), p. 115767. DOI: https://doi.org/10.1016/j.cam.2024.115767.

[6] D. J. Knezevic and J. W. Peterson. "A high-performance parallel implementation of the certified reduced basis method". In: *Computer Methods in Applied Mechanics and Engineering* 200.13 (2011), pp. 1455–1466. DOI: https://doi.org/10.1016/j.cma.2010.12.026.

[7] M. Drohmann, B. Haasdonk, S. Kaulmann, and M. Ohlberger. "A Software Framework for Reduced Basis Methods Using Dune-RB and RBmatlab". In: Jan. 2012. DOI: 10.1007/978-3-642-28589-9_6.

[8] Y. Choi, W. J. Arrighi, D. M. Copeland, R. W. Anderson, and G. M. Oxberry. *libROM*. Lawrence Livermore National Laboratory.

[9] F. Ballarin, A. Sartori, and G. Rozza. "RBniCS - reduced order modelling in FEniCS". In: *ScienceOpen Posters* (2015).

[10] R. Fritze. "pyMOR - Model Order Reduction with Python". In: July 2014. DOI: 10.6084/m9.figshare.1134465.

[11] Y. Kim, K. Wang, and Y. Choi. "Efficient Space–Time Reduced Order Model for Linear Dynamical Systems in Python Using Less than 120 Lines of Code". In: *Mathematics* 9.14 (July 2021), p. 1690. DOI: 10.3390/math9141690.

[12] A. Logg and G. N. Wells. "DOLFIN: Automated finite element computing". In: *ACM Transactions on Mathematical Software* 37.2 (Apr. 2010), pp. 1–28. DOI: 10.1145/1731022.1731030.

[13] D. Arndt et al. "The deal.II finite element library: Design, features, and insights". In: *Computers & Mathematics with Applications* 81 (2021). Development and Application of Open-source Software for Problems with Numerical PDEs, pp. 407–422. DOI: https://doi.org/10.1016/j.camwa.2020.02.022.

[14] J. Bezanson, A. Edelman, S. Karpinski, and V. B. Shah. "Julia: A Fresh Approach to Numerical Computing". In: *SIAM Review* 59.1 (2017), pp. 65–98. DOI: 10.1137/141000671. eprint: https://doi.org/10.1137/141000671.

[15] S. Badia and F. Verdugo. "Gridap: An extensible Finite Element toolbox in Julia". In: *Journal of Open Source Software* 5.52 (2020), p. 2520. DOI: 10.21105/joss.02520.

[16] F. Verdugo and S. Badia. "The software design of Gridap: a Finite Element package based on the Julia JIT compiler". In: *Computer Physics Communications* 276 (2022), p. 108341. DOI: 10.1016/j.cpc.2022.108341.

[17] M. Lubin et al. "JuMP 1.0: Recent improvements to a modeling language for mathematical optimization". In: *Mathematical Programming Computation* (2023). DOI: 10.1007/s12532-023-00239-3.

[18] C. Rackauckas and Q. Nie. "DifferentialEquations.jl–a performant and feature-rich ecosystem for solving differential equations in Julia". In: *Journal of Open Research Software* 5.1 (2017).

[19] F. Negri, A. Manzoni, and G. Rozza. "Reduced basis approximation of parametrized optimal flow control problems for the Stokes equations". In: *Computers & Mathematics with Applications* 69.4 (2015), pp. 319–336.

[20] N. Halko, P.-G. Martinsson, and J. A. Tropp. "Finding structure with randomness: probabilistic algorithms for constructing approximate matrix decompositions". In: *SIAM review* 53.2 (2011), pp. 217–288.

[21] K. Ho et al. *JuliaMatrices/LowRankApprox.jl: v0.5.2*. Version v0.5.2. Mar. 2022. DOI: 10.5281/zenodo.6394438.

[22] D. Amsallem, M. J. Zahr, and C. Farhat. "Nonlinear model order reduction based on local reduced-order bases". In: *International Journal for Numerical Methods in Engineering* 92.10 (2012), pp. 891–916. DOI: https://doi.org/10.1002/nme.4371. eprint: https://onlinelibrary.wiley.com/doi/pdf/10.1002/nme.4371.

[23] B. Peherstorfer, D. Butnaru, K. Willcox, and H.-J. Bungartz. "Localized Discrete Empirical Interpolation Method". In: *SIAM Journal on Scientific Computing* 36.1 (2014), A168–A192. DOI: 10.1137/130924408. eprint: https://doi.org/10.1137/130924408.

[24] S. Pagani, A. Manzoni, and A. Quarteroni. "Numerical approximation of parametrized problems in cardiac electrophysiology by a local reduced basis method". In: *Computer Methods in Applied Mechanics and Engineering* 340 (2018), pp. 530–558. DOI: https://doi.org/10.1016/j.cma.2018.06.003.

[25] Y. Choi, G. Boncoraglio, S. Anderson, D. Amsallem, and C. Farhat. *Gradient-based Constrained Optimization Using a Database of Linear Reduced-Order Models*. 2020. arXiv: 1506.07849 [math.NA].

[26] M. Chasapi, P. Antolin, and A. Buffa. "A localized reduced basis approach for unfitted domain methods on parameterized geometries". In: *Computer Methods in Applied Mechanics and Engineering* 410 (2023), p. 115997. DOI: https://doi.org/10.1016/j.cma.2023.115997.

[27] N. Mueller, Y. Zhao, S. Badia, and T. Cui. "A tensor-train reduced basis solver for parameterized partial differential equations on Cartesian grids". In: *Journal of Computational and Applied Mathematics* 472 (2026), p. 116790. DOI: https://doi.org/10.1016/j.cam.2025.116790.

[28] J. T. Lauzon et al. "S-OPT: A Points Selection Algorithm for Hyper-Reduction in Reduced Order Models". In: *SIAM Journal on Scientific Computing* 46.4 (July 2024), B474–B501. DOI: 10.1137/22m1484018.

[29] S. Salsa. *Partial Differential Equations in Action*. Vol. 99. Springer, 2016.

[30] A. Quarteroni. *Numerical Models for Differential Problems*. Vol. 8. Springer, 2016.

[31] Y. Choi and K. Carlberg. "Space–time least-squares Petrov–Galerkin projection for nonlinear model reduction". In: *SIAM Journal on Scientific Computing* 41.1 (2019), A26–A58.

[32] I. Oseledets and E. Tyrtyshnikov. "TT-cross approximation for multidimensional arrays". In: *Linear Algebra and its Applications* 432.1 (2010), pp. 70–88.

[33] I. V. Oseledets. "Tensor-train decomposition". In: *SIAM Journal on Scientific Computing* 33.5 (2011), pp. 2295–2317.

[34] S. Chaturantabut and D. C. Sorensen. "Nonlinear model reduction via discrete empirical interpolation". In: *SIAM Journal on Scientific Computing* 32.5 (2010), pp. 2737–2764.

[35] 2. Weiwu Jiang Xiaowei Gao. "Review of Collocation Methods and Applications in Solving Science and Engineering Problems". In: *Computer Modeling in Engineering & Sciences* 140.1 (2024), pp. 41–76. DOI: 10.32604/cmes.2024.048313.

[36] Y. Choi, P. Brown, W. Arrighi, R. Anderson, and K. Huynh. "Space–time reduced order model for large-scale linear dynamical systems with application to Boltzmann transport problems". In: *Journal of Computational Physics* 424 (2021).

[37] M. Pharr, W. Jakob, and G. Humphreys. "07 - Sampling and Reconstruction". In: *Physically Based Rendering (Third Edition)*. Ed. by M. Pharr, W. Jakob, and G. Humphreys. Third Edition. Boston: Morgan Kaufmann, 2017, pp. 401–504. DOI: https://doi.org/10.1016/B978-0-12-800645-0.50007-5.

[38] M. D. McKay, R. J. Beckman, and W. J. Conover. "A Comparison of Three Methods for Selecting Values of Input Variables in the Analysis of Output from a Computer Code". In: *Technometrics* 21.2 (1979), pp. 239–245.

[39] D. Boffi, F. Brezzi, M. Fortin, et al. *Mixed finite element methods and applications*. Vol. 44. Springer, 2013.

[40] G. Rozza. "On optimization, control and shape design of an arterial bypass". In: *International Journal for Numerical Methods in Fluids* 47.10-11 (2005), pp. 1411–1419.

[41] F. Ballarin, A. Manzoni, A. Quarteroni, and G. Rozza. "Supremizer stabilization of POD–Galerkin approximation of parametrized steady incompressible Navier–Stokes equations". In: *International Journal for Numerical Methods in Engineering* 102.5 (2015), pp. 1136–1161.

[42] N. Dal Santo and A. Manzoni. "Hyper-reduced order models for parametrized unsteady Navier–Stokes equations on domains with variable shape". In: *Advances in Computational Mathematics* 45.5 (2019), pp. 2463–2501.

[43] L. Pegolotti, M. R. Pfaller, A. L. Marsden, and S. Deparis. "Model order reduction of flow based on a modular geometrical approximation of blood vessels". In: *Computer methods in applied mechanics and engineering* 380 (2021), p. 113762.

[44] N. Mueller, S. Badia, and Y. Zhao. *Reduced basis solvers for unfitted methods on parameterized domains*. 2025. arXiv: 2508.15320 [math.NA].

[45] A. Hansbo and P. Hansbo. "An unfitted finite element method, based on Nitsche's method, for elliptic interface problems". In: *Computer Methods in Applied Mechanics and Engineering* 191.47 (2002), pp. 5537–5552. DOI: https://doi.org/10.1016/S0045-7825(02)00524-8.

[46] E. Burman, S. Claus, P. Hansbo, M. G. Larson, and A. Massing. "CutFEM: Discretizing geometry and partial differential equations". In: *International Journal for Numerical Methods in Engineering* 104.7 (2015), pp. 472–501. DOI: https://doi.org/10.1002/nme.4823. eprint: https://onlinelibrary.wiley.com/doi/pdf/10.1002/nme.4823.

[47] D. Elfverson, M. G. Larson, and K. Larsson. "CutIGA with basis function removal". In: *Advanced Modeling and Simulation in Engineering Sciences* 5.1 (Mar. 2018). DOI: 10.1186/s40323-018-0099-2.

[48] T. Belytschko, N. Moës, S. Usui, and C. Parimi. "Arbitrary discontinuities in finite elements". In: *International Journal for Numerical Methods in Engineering* 50.4 (2001), pp. 993–1013. DOI: https://doi.org/10.1002/1097-0207(20010210)50:4<993::AID-NME164>3.0.CO;2-M. eprint: https://onlinelibrary.wiley.com/doi/pdf/10.1002/1097-0207%2820010210%2950%3A4%3C993%3A%3AAID-NME164%3E3.0.CO%3B2-M.

[49] S. Badia, F. Verdugo, and A. F. Martín. "The aggregated unfitted finite element method for elliptic problems". In: *Computer Methods in Applied Mechanics and Engineering* 336 (2018), pp. 533–553. DOI: https://doi.org/10.1016/j.cma.2018.03.022.

[50] S. Badia, A. F. Martín, and F. Verdugo. “GridapDistributed: a massively parallel finite element toolbox in Julia”. In: *Journal of Open Source Software* 7.74 (2022), p. 4157. DOI: 10.21105/joss.04157.
[51] F. Verdugo et al. *PartitionedArrays/PartitionedArrays.jl: v0.3.5*. Version v0.3.5. Aug. 2025. DOI: 10.5281/zenodo.16793213.
[52] S. Balay et al. *PETSc Web page*. http://www.mcs.anl.gov/petsc. 2015.
[53] J. Manyer, A. F. Martín, and S. Badia. “GridapSolvers.jl: Scalable multiphysics finite element solvers in Julia”. In: *Journal of Open Source Software* 9.102 (2024), p. 7162. DOI: 10.21105/joss.07162.
[54] Y. Kim, Y. Choi, D. Widemann, and T. Zohdi. “A fast and accurate physics-informed neural network reduced order model with shallow masked autoencoder”. In: *Journal of Computational Physics* 451 (Feb. 2022), p. 110841. DOI: 10.1016/j.jcp.2021.110841.
[55] A. N. Diaz, Y. Choi, and M. Heinkenschloss. “A fast and accurate domain decomposition nonlinear manifold reduced order model”. In: *Computer Methods in Applied Mechanics and Engineering* 425 (May 2024), p. 116943. DOI: 10.1016/j.cma.2024.116943.
[56] K. Stein, T. E. Tezduyar, and R. J. Benney. “Mesh Moving Techniques for Fluid-Structure Interactions with Large Displacements”. In: *Journal of Applied Mechanics* 70.1 (2003), pp. 58–63. DOI: 10.1115/1.1530635.
[57] K. Stein, T. E. Tezduyar, and R. J. Benney. “Automatic Mesh Update with the Solid-Extension Mesh Moving Technique”. In: *Computer Methods in Applied Mechanics and Engineering* 193.21 (2004), pp. 2019–2032. DOI: 10.1016/j.cma.2004.01.023.
[58] T. E. Tezduyar, M. Behr, S. Mittal, and A. A. Johnson. “Computation of Unsteady Incompressible Flows with the Stabilized Finite Element Methods: Space-Time Formulations, Iterative Strategies and Massively Parallel Implementations”. In: *New Methods in Transient Analysis, AMD-Vol. 246*. ASME, 1992, pp. 7–24.
[59] A. Nouy, M. Chevreuil, and E. Safatly. “Fictitious domain method and separated representations for the solution of boundary value problems on uncertain parameterized domains”. In: *Computer Methods in Applied Mechanics and Engineering* 200.45 (2011), pp. 3066–3082. DOI: https://doi.org/10.1016/j.cma.2011.07.002.
[60] M. Balajewicz and C. Farhat. “Reduction of nonlinear embedded boundary models for problems with evolving interfaces”. In: *Journal of Computational Physics* 274 (2014), pp. 489–504. DOI: https://doi.org/10.1016/j.jcp.2014.06.038.
[61] E. N. Karatzas, F. Ballarin, and G. Rozza. “Projection-based reduced order models for a cut finite element method in parametrized domains”. In: *Computers & Mathematics with Applications* 79.3 (2020), pp. 833–851. DOI: https://doi.org/10.1016/j.camwa.2019.08.003.
[62] E. N. Karatzas, M. Nonino, F. Ballarin, and G. Rozza. “A Reduced Order Cut Finite Element method for geometrically parametrized steady and unsteady Navier–Stokes problems”. In: *Computers & Mathematics with Applications* 116 (2022). New trends in Computational Methods for PDEs, pp. 140–160. DOI: https://doi.org/10.1016/j.camwa.2021.07.016.
[63] F. de Prenter, C. Verhoosel, G. van Zwieten, and E. van Brummelen. “Condition number analysis and preconditioning of the finite cell method”. In: *Computer Methods in Applied Mechanics and Engineering* 316 (2017). Special Issue on Isogeometric Analysis: Progress and Challenges, pp. 297–327. DOI: https://doi.org/10.1016/j.cma.2016.07.006.

## Appendix A. Code used to solve the unsteady Navier Stokes application

In this section, we present the code used to generate the results in Section 4. To improve readability, we separate the definition of FE-related quantities from the variables used in the ROM routines. These two components are shown in Listings 13 and 14, respectively.

We first discuss the FE code. The geometry in Fig. 11 is imported from file, and physical labels are introduced (lines 10-19). These labels are used to manage boundary conditions appropriately: in particular, we define the no-slip, no-penetration, and inlet Dirichlet boundary conditions detailed in Section 4 on the correct physical groups (lines 78-84). At lines 91-94, we define multi-field test and trial spaces using the `BlockMultiFieldStyle` identifier, which constructs FE matrices and vectors in a block-wise layout rather than storing their entries consecutively. This structure is required in GridapROMs for efficiency, as reduced subspaces are constructed separately for each unknown variable (velocity and pressure in this example), and the block-wise layout naturally facilitates this separation. The nonlinear convective term is separated from the linear forms, with two distinct FE operators defined (lines 97-100) and combined into a single `LinearNonlinearTransientParamOperator` (lines 101-102). This separation serves two purposes: first, it optimizes the time-marching scheme (the Crank-Nicolson scheme defined at lines 21-26, 111-113) by avoiding repeated assembly of FE matrices and vectors associated with the linear forms at each Newton-Raphson iteration for a fixed time step; second, empirical evidence shows that treating linear and nonlinear components separately improves hyper-reduction approximations. The remaining lines of Listing 13 follow the implementation of a standard GridapROMs driver and are omitted here to avoid repetition.

Next, we briefly discuss the implementation of the ROM driver shown in Listing 14. We provide the `RBSolver` with the construction of a $(H^1(\Omega))^3$-orthogonal subspace for the velocity, and a $L^2(\Omega)$-orthogonal subspace for the pressure. We recall that, by passing `coupling` to the solver implicitly requests the enrichment of the velocity subspace

lst_navier_stokes_fe.jl

```julia
# packages
using GridapROMs
using Gridap
using Gridap.MultiField
using GridapSolvers
using GridapSolvers.LinearSolvers
using GridapSolvers.NonlinearSolvers
using DrWatson

# load geometry
model_dir = datadir(joinpath("models","circle.json"))
Ωₕ = DiscreteModelFromFile(model_dir)
τₕ = Triangulation(Ωₕ)

# set boundary labels
labels = get_face_labeling(Ωₕ)
add_tag_from_tags!(labels,"noslip",["cylinders","walls"])
add_tag_from_tags!(labels,"nopenetration",["top_bottom"])
add_tag_from_tags!(labels,"dirichlet",["inlet"])

# time marching
θ = 0.5
dt = 0.0025
t₀ = 0.0
tf = 60dt
tdomain = t₀:dt:tf

# parametric quantities
pdomain = (1,10,1,10,1,10)
D = TransientParamSpace(pdomain,tdomain)

# viscosity
const Re = 100.0
ν(μ,t) = x -> μ[1]/Re*exp((sin(t)+cos(t))/sum(μ))
νₚₜ(μ,t) = parameterize(ν,μ,t)

# Dirichlet data
const W = 0.5
flow(μ,t) = abs(1-cos(π*t/T)+μ[3]sin(μ[2]π*t/T)/100)
g(μ,t) = x -> VectorValue(x[2]*(x[2]-W)*flow(μ,t),0.0,0.0)
gₚₜ(μ,t) = parameterize(g,μ,t)
g₀(μ,t) = x -> VectorValue(0.0,0.0,0.0)
g₀ₜ(μ,t) = parameterize(g₀,μ,t)

# numerical integration
order = 2
dΩₕ = Measure(τₕ,2order+1)

# weak form - linear terms
m(μ,t,(uₜ,pₜ),(v,q),dΩₕ) = ∫(v⋅uₜ)dΩₕ
a(μ,t,(u,p),(v,q),dΩₕ) = (
    ∫(νₚₜ(μ,t)*∇(v)⊙∇(u))dΩₕ -
    ∫(p*(∇⋅(v)))dΩₕ +
    ∫(q*(∇⋅(u)))dΩₕ
)
res_lin(μ,t,(u,p),(v,q),dΩₕ) = (
    m(μ,t,(∂t(u),∂t(p)),(v,q),dΩₕ) +
    a(μ,t,(u,p),(v,q),dΩₕ)
)

# weak form - nonlinear terms
conv(u,∇u) = (∇u')⋅u
res_nlin(μ,t,(u,p),(v,q),dΩₕ) = (
    ∫( v⊙(conv∘(u,∇(u))) )dΩₕ
)
dconv(du,∇du,u,∇u) = conv(u,∇du)+conv(du,∇u)
jac_nlin(μ,t,(u,p),(du,dp),(v,q),dΩₕ) = (
    ∫( v⊙(dconv∘(du,∇(du),u,∇(u))) )dΩₕ
)

# triangulation information
τₕ_jac = (τₕ,)
τₕ_djac = (τₕ,)
τₕ_res = (τₕ,)

# FE interpolation - velocity
reffeᵤ = ReferenceFE(lagrangian,VectorValue{3,Float64},order)
V = TestFESpace(Ωₕ,reffeᵤ;conformity=:H1,
    dirichlet_tags=["dirichlet","nopenetration","noslip"],
    dirichlet_masks=[
    [true,true,true],
    [false,false,true],
    [true,true,true]])
U = TransientTrialParamFESpace(V,[gₚₜ,gₚₜ,g₀ₜ])

# FE interpolation - pressure
reffeₚ = ReferenceFE(lagrangian,Float64,order)
Q = TestFESpace(Ωₕ,reffeₚ;conformity=:C0)
P = TrialFESpace(Q)

# FE interpolation - multi field
style = BlockMultiFieldStyle()
Y = TransientMultiFieldParamFESpace([V,Q];style)
X = TransientMultiFieldParamFESpace([U,P];style)

# linear - nonlinear splitting
feop_lin = TransientParamLinearFEOperator(res_lin,a,m,
    D,X,Y,τₕ_res,(τₕ_jac,τₕ_djac))
feop_nlin = TransientParamFEOperator(res_nlin,jac_nlin,
    D,X,Y,τₕ_res,(τₕ_jac,))
feop = LinearNonlinearTransientParamFEOperator(
    feop_lin,feop_nlin)

# initial condition
u₀(x,μ) = x -> VectorValue(0.0,0.0,0.0)
u₀ₚ(μ) = parameterize(u₀,μ)
p₀(x,μ) = x -> 0.0
p₀ₚ(μ) = parameterize(p₀,μ)
xh₀ₚ(μ) = interpolate_everywhere([u₀ₚ(μ),p₀ₚ(μ)],X(μ,t₀))

# FE solver
nlslvr = NewtonSolver(LUSolver();rtol=1e-10,maxiter=20)
odeslvr = ThetaMethod(nlslvr,dt,θ)
```

FIGURE 13. FE code used to solve the unsteady Navier Stokes application with GridapROMs.

with supremizers. The keyword `sketch` is used to select a randomized POD technique for compressing the snapshots, implemented via the Julia package LowRankApprox [21]; in this case, the chosen randomization method corresponds to a Sparse Random Gaussian stencil [20]. The remaining code closely follows the usage example presented in Subsection 2.5.

## APPENDIX B. CODE USED TO SOLVE A POISSON EQUATION ON A PARATERIZED DOMAIN USING A LOCAL TT-RB APPROACH

In this section, we present the code used to solve a Poisson equation on a parameterized domain using a local TT-based approach. The application consists of a square with a circular hole of fixed radius, free to move along a portion of the

```
lst_navier_stokes_rb.jl

# RB solver
tol = 1e-4
coupling((u,p),(v,q)) = ∫(p*(∇⋅(v)))dΩₕ
inner_prod((du,dp),(v,q)) = (
    ∫(u⋅v)dΩₕ + ∫(∇(v)⊙∇(u))dΩₕ + ∫(p*q)dΩₕ
)
red_sol = TransientReduction(coupling,tol,inner_prod;
    nparams=55,sketch=:sprn)
rbslvr = RBSolver(odeslvr,red_sol;
    nparams_res=55,nparams_jac=15,nparams_djac=1)

dir = datadir("navier_stokes")
create_dir(dir)

rbop = try
    # try loading offline quantities
    load_operator(dir,feop)
catch
    # compute and save offline quantities
    reduced_operator(dir,rbslvr,feop,xh₀ₚ)
end

# online phase
μₒₙ = realization(feop;nparams=5,sampling=:uniform)
x̂,stats = solution_snapshots(rbslvr,feop,μₒₙ,xh₀ₚ)

# post process
x,rbstats = solve(odeslvr,rbop,μₒₙ,xh₀ₚ)
rb_performance(rbslvr,feop,rbop,x̂,x,stats,rbstats)
```

FIGURE 14. ROM code used to solve the unsteady Navier Stokes application with GridapROMs.

square's main diagonal. The hole's displacement along the diagonal is described by a single shape parameter. This setup is nontrivial, both conceptually and in terms of implementation. The strategy closely follows the methodology introduced in [44], to which we refer the interested reader for further details. We briefly summarize the procedure below:

- We use mesh deformation maps from a reference configuration – where the circle is centered at the midpoint of the diagonal – to each physical configuration. This allows us to define a single pair of test-trial FE spaces on the reference domain, and to integrate the weak form of the PDE via pull-back mappings. In this example, elastic deformations [56–58] are employed to map the reference configuration onto the physical ones.
- At the HF level, we rely on an unfitted FE discretization [45–49], in which the circle's boundary is implicitly represented through level-sets defined on a Cartesian background grid. This boundary is thus referred to as *embedded*. Compared to the drivers in Listings 6 and 13, two key differences arise. First, when the embedded boundary is subject to Dirichlet conditions, the data must be imposed weakly, e.g. using Nitsche's penalty method. Second, the geometry and FE spaces require specialized types, handled by the GridapEmbedded package.
- The solution of the PDE is defined only on the active portion of the background mesh. To apply TT-RB, the solution must be extended to the remaining *external* cells, ensuring that snapshots can be represented as a 3D tensor (two spatial dimensions plus one parametric axis). Extension strategies for unfitted discretizations have been widely studied [59–62]. A simple approach is the zero extension, but higher-quality reduced subspaces are obtained with harmonic extensions, which we adopt in this example.
- Finally, we employ a local TT-RB strategy for the reduced approximation of the Poisson equation. This strategy is selected through the commands described in Subsections 3.2-3.3. A tighter tolerance is used in the hyper-reduction of both RHS and LHS, so that the accuracy is primarily governed by the projection error.

Having outlined the general aspects of the benchmark, we now turn to the specifics of the code. For clarity, the implementation is split into three listings: (i) Listing 15, where we define the unfitted reference configuration and the elastic deformations; (ii) Listing 16, where we implement the pull-back mappings for integration over the reference configuration and set up the interfaces required for the harmonic extension of the solutions; and (iii) Listing 17, where the ROM routines are executed.

We now proceed step by step. In Listing 15, we begin by defining the (background) discrete model $\Omega_h$_bg over the square $[0,2]^2$ (line 14), using a `TProductDiscreteModel` (see Subsection 3.2 for more details). Parametric forcing and Dirichlet terms are then introduced in the usual way.

The new elements appear starting at line 24, where quantities related to the unfitted reference configuration are set up. Between lines 24-30, we cut a disk of radius 0.3 centered at $(1,1)$ from $\Omega_h$_bg. Subsequently, we define the background, active, and physical triangulations, the cut-interface triangulation and its normal vector, together with their associated measures (lines 32-43). Lines 45-47 introduce the cell aggregation method described in [49], which is used throughout the benchmark to construct the unfitted FE spaces. In principle, only small cut cells – those that adversely affect the conditioning of the algebraic systems arising from unfitted methods (see [63]) – need to be aggregated. Here, the identifier `AggregateAllCutCells` aggregates all cut cells, though in practice users may prefer to aggregate only those below a given size threshold. The function `get_deformation_map` sets up the linear elasticity problem used to compute the deformation map. Specifically, the map is obtained by prescribing a known boundary displacement $\underline{\varphi}^{\mu}(\underline{x}) \doteq (\mu_1 - 1, \mu_1 - 1)$ weakly on $\Gamma_h$ – the circle's boundary – while enforcing zero displacement strongly on the outer boundaries of the square. Note that the elastic problem is solved with a higher interpolation order than the Poisson equation. Finally, at line 93 we encounter the first instance of an (aggregated) unfitted FE space, constructed via the function `AgFEMSpace`. This function takes as input a standard FE space defined on the active triangulation along with the cell aggregates introduced earlier. The

```
moving_domains_deformation.jl

using DrWatson
using Gridap
using GridapEmbedded
using GridapROMs

# geometry
const L = 2
const W = 2
const N = 40
const h = max(L,W)/N

domain = (0,L,0,W)
partition = (N,N)
Ωₕ_bg = TProductDiscreteModel(domain,partition)

# parametric quantities
pdomain = (0.6,1.4)
D = ParamSpace(pdomain)
f(μ) = x -> 1.0
g(μ) = x -> 0.0
fₚ(μ) = parameterize(f,μ)
gₚ(μ) = parameterize(g,μ)

# unfitted geometry: circular hole cut from the
# background domain; we start by defining the
# reference configuration
μ₀ = (1.0,)
x₀ = Point(μ₀[1],μ₀[1])
geo = !disk(0.3,x0=x₀)
cutgeo = cut(Ωₕ_bg,geo)

τₕ_bg = Triangulation(Ωₕ_bg)
τₕ_act = Triangulation(cutgeo,ACTIVE)
τₕ = Triangulation(cutgeo,PHYSICAL)
Γₕ = EmbeddedBoundary(cutgeo)
n = get_normal_vector(Γₕ)

# numerical integration
order = 1
degree = 2*order
dΩₕ_bg = Measure(τₕ_bg,degree)
dΩₕ = Measure(τₕ,degree)
dΓₕ = Measure(Γₕ,degree)

# Aggregated cell method
strategy = AggregateAllCutCells()
aggregates = aggregate(strategy,cutgeo)

# Nitsche penalty method
const γ = 10.0
const Cₕ = γ/h

# FE solver
slvr = LUSolver()

# linear elasticity deformation to use for the
# deformed configuration
function get_deformation_map(μ)
    # translation of the circular hole
    φ(μ) = x -> VectorValue(μ[1]-μ₀[1],μ[1]-μ₀[1])
    φₚ(μ) = parameterize(φ,μ)

    # linear elasticity coefficients
    E = 1
    ν = 0.33
    λ = E*ν/((1+ν)*(1-2*ν))
    p = E/(2(1+ν))

    # linear constitutive law
    σ(ε) = λ*tr(ε)*one(ε) + 2*p*ε

    # set a higher order for the linear elasticity
    order = 2
    degree = 2*order
    dΩₕ = Measure(τₕ,degree)
    dΓₕ = Measure(Γₕ,degree)

    # weak form
    a(μ,u,v) = (
        ∫( ε(v) ⊙ (σ∘ε(u)) )*dΩₕ +
        ∫( Cₕ*v⋅u - v⋅(n⋅(σ∘ε(u))) - (n⋅(σ∘ε(v)))⋅u )dΓₕ
    )
    l(μ,v) = ∫( Cₕ*v⋅φₚ(μ) - (n⋅(σ∘ε(v)))⋅φₚ(μ) )dΓₕ
    res(μ,u,v) = ∫( ε(v) ⊙ (σ∘ε(u)) )*dΩₕ - l(μ,v)

    # FE interpolation
    reffeφ = ReferenceFE(
        lagrangian,VectorValue{2,Float64},order
    )
    Vφ_act = FESpace(τₕ_act,reffeφ,
        conformity=:H1,dirichlet_tags="boundary"
    )
    Vφ_agg = AgFEMSpace(Vφ_act,aggregates)
    Uφ_agg = ParamTrialFESpace(Vφ_agg)

    # solve
    feop = LinearParamFEOperator(res,a,D,Uφ_agg,Vφ_agg)
    solve(LinearFESolver(slvr),feop,μ)
end
```

FIGURE 15. ROM code used to solve the unsteady Navier Stokes application with GridapROMs.

parameter $\boldsymbol{\mu}$ is only used at the last line of the function to actually run the solver.

Next, we turn to the code in Listing 16. Most of the implementation resides in `def_extended_fe_operator`, a function which builds an FE operator supporting the harmonic extension to the exterior cells for a given parameter $\boldsymbol{\mu}$. The procedure begins with the computation of the deformation map via a call to `get_deformation_map`, followed by the construction of a `CellField` from the cell map associated with the deformed triangulation (lines 28-35). This `CellField` is then employed, between lines 37-58, to define the following pull-backs:

$$\boldsymbol{\nabla}^{\mu} \doteq \left(\underline{\underline{J}}^{\mu}\right)^{-T} \boldsymbol{\nabla}, \quad \int_{\Omega^{\mu}} d\Omega^{\mu} \doteq \int_{\Omega} \det(\underline{\underline{J}}^{\mu}) d\Omega, \quad \underline{n}^{\mu} \doteq \frac{\left(\underline{\underline{J}}^{\mu}\right)^{-T} \underline{n}}{\| \left(\underline{\underline{J}}^{\mu}\right)^{-T} \underline{n} \|_2}, \quad \int_{\Gamma^{\mu}} d\Gamma^{\mu} \doteq \int_{\Gamma} \| \left(\underline{\underline{J}}^{\mu}\right)^{-T} \underline{n} \|_2 \det(\underline{\underline{J}}^{\mu}) d\Gamma, \tag{17}$$

where:

- $\underline{\underline{J}}^{\mu} \doteq \boldsymbol{\nabla}^{\mu} \underline{\varphi}^{\mu} \in \mathbb{R}^{2\times 2}$ is the Jacobian of the deformation map;

moving_domains_fe.jl

```julia
using Gridap.CellData
using Gridap.Geometry
import Gridap.Geometry: push_normal

using GridapROMs.Extensions

# triangulation information
τₕ_l = (τₕ,)
τₕ_a = (τₕ,Γₕ)
domains = FEDomains(τₕ_l,τₕ_a)

# FE interpolation
reffe = ReferenceFE(lagrangian,Float64,order)
V_act = FESpace(τₕ_act,reffe,
    conformity=:H1,dirichlet_tags=[1,3,7]
)
V_agg = AgFEMSpace(V_act,aggregates)

# extension FE solver: we solve the actual PDE
# with the input solver slvr, and additionally
# we extend the nodal values to cover the entire
# background mesh (i.e., inside the circular hole).
# By default, a harmonic extension is chosen. This
# will allow us to run TT-RB
slvr_ext = ExtensionSolver(slvr)

function def_extended_fe_operator(μ)
    # linear elasticity deformation
    φₕ = get_deformation_map(μ)
    # deformed triangulation
    τφₕ_act = mapped_grid(τₕ_act,φₕ)

    # cell map of the deformed triangulation
    ϕ = get_cell_map(τφₕ_act)
    ϕₕ = GenericCellField(ϕ,τₕ_act,ReferenceDomain())

    # determinant of the deformation's Jacobian
    dJφ = det∘(∇(ϕₕ))

    # inverse of the deformation's Jacobian
    ϕₕ⁻¹ = inv∘∇(ϕₕ)

    # gradient in the deformed configuration
    ∇φ(a) = dot∘(ϕₕ⁻¹,∇(a))

    # normal vector to the circular hole in the
    # deformed configuration
    nφ = push_normal∘(ϕₕ⁻¹,n)

    # determinant of the deformation's Jacobian
    # restricted to the boundary of the circular hole
    # in the deformed configuration
    ϕₕϕₕᵀ = (j->j⋅j')∘(∇(ϕₕ))
    dJφ_Γ = ((j,c,n)->j*√(n⋅inv(c)⋅n))∘(dJφ,ϕₕϕₕᵀ,nφ)

    # integrals in the deformed configuration
    ∫_Ω(a) = ∫(a*dJφ)
    ∫_Γ(a) = ∫(a*dJφ_Γ)

    # PDE in the deformed configuration
    a(μ,u,v,dΩₕ,dΓₕ) = (
        ∫_Ω(∇φ(v)⋅∇φ(u))dΩₕ +
        ∫_Γ( Cₕ*v*u - v*(nφ⋅∇φ(u)) - (nφ⋅∇φ(v))*u )dΓₕ
    )
    l(μ,v,dΩₕ) = ∫_Ω(fₚ(μ)⋅v)dΩₕ
    res(μ,u,v,dΩₕ) = ∫_Ω(∇φ(v)⋅∇φ(u))dΩₕ - l(μ,v,dΩₕ)

    # FE interpolation in the deformed configuration
    V_bg = FESpace(τₕ_bg,reffe,
        conformity=:H1,dirichlet_tags=[1,3,7]
    )
    # the following FE spaces include information related to
    # circular hole, allowing for the extension of the nodal
    # values to take place. This will allow us to run TT-RB
    V_ext = DirectSumFESpace(V_bg,V_agg)
    U_ext = ParamTrialFESpace(V_ext,gₚ)
    ExtensionLinearParamOperator(res,a,D,U_ext,V_ext,domains)
end
```

FIGURE 16. ROM code used to solve the unsteady Navier Stokes application with GridapROMs.

- $\boldsymbol{\nabla}^{\mu}$ (resp. $\boldsymbol{\nabla}$) is the gradient operator on the deformed (resp. reference) domain;
- $d\Omega^{\mu}$ (resp. $d\Omega$) is the volume measure on the deformed (resp. reference) domain;
- $d\Gamma^{\mu}$ (resp. $d\Gamma$) is the surface measure on the deformed (resp. reference) boundary;
- $\underline{n}^{\mu}$ (resp. $\underline{n}$) is the outward unit normal vector on the deformed (resp. reference) boundary.

The pull-backs in (17) are then employed to formulate the weak form of the Poisson problem at lines 60-69. We plan to fully automate the implementation of the pull-backs in future releases of GridapROMs, so that users will no longer need to define them manually. The remaining novelties of this script concern the definition of a nonstandard pair of FE spaces and a numerical solver that support the harmonic extension to the exterior cells. The construction of the solver is straightforward: one simply calls the constructor `ExtensionSolver`, which takes the previously defined `LUSolver` as input. The definition of the spaces, however, is more involved and requires two components:

- A FE space (not necessarily aggregate) for the Poisson problem (lines 12-17).
- A (background) FE space defined on the entire background domain (lines 68-71).

Note that the former alone would suffice if the extension procedure were not required. Since it is, the background FE space – which provides the necessary information for the exterior cells – must also be provided. With these ingredients, we construct a pair of `DirectSumFESpaces`, which combine the background FE space and the aggregate FE space while internally building a “complementary” FE space on the exterior cells. The “direct sum” of the aggregate and complementary spaces recovers the background space. The complementary space is used for the harmonic extension of solutions defined on the aggregate space. Finally, a `ExtensionLinearParamOperator` is returned, taking this pair of spaces as input.

Lastly, we illustrate how the ROM routines are executed for this application. At lines 3-14, we set up a local TT-RB solver with tolerances of $10^{-4}$ in both spatial directions, a randomized compression scheme, and $8$ centroids for the `Kmeans` clustering of the 100 offline parameters. The same settings are applied to the hyper-reduction step, except for the

```
moving_domains_rb.jl

 1 using GridapROMs.RBSteady
 2
 3 # local TT-RB solver
 4 tol = fill(1e-4,2)
 5 nparams = 100
 6 ncentroids = 8
 7 kwargs = (
 8     nparams=nparams,
 9     sketch=:sprn,
10     ncentroids=ncentroids
11 )
12
13 inner_prod(du,v) = ∫(∇(v)⋅∇(du))dΩₕ_bg
14 red_sol = LocalReduction(tol,inner_prod;kwargs...)
15 red_res = LocalHyperReduction(tol.*1e-2;kwargs...)
16 red_jac = LocalHyperReduction(tol.*1e-2;kwargs...)
17 rbslvr = RBSolver(slvr_ext,red_sol,red_res,red_jac)
18
19 # offline phase
20 μ = realization(D;nparams)
21 feop = def_extended_fe_operator(μ)
22 fesnaps, = solution_snapshots(rbslvr,feop,μ)
23 rbop = reduced_operator(rbslvr,feop,fesnaps)
24
25 # solve to use when using the localization strategy
26 # on a moving domain
27 function local_rb_solve(rbslvr,rbop,μ)
28     trial = get_trial(rbop)
29     k, = get_clusters(trial)
30     μsplit = cluster(μ,k)
31     perfs = ROMPerformance[]
32     for μᵢ in μsplit
33         feopᵢ = def_extended_fe_operator(μᵢ)
34         xᵢ,festatsᵢ = solution_snapshots(rbslvr,feopᵢ,μᵢ)
35         rbopᵢ = change_operator(
36             get_local(rbop,first(μᵢ)),feopᵢ
37         )
38         x̂ᵢ,rbstatsᵢ = solve(rbslvr,rbopᵢ,μᵢ)
39         perfᵢ = eval_performance(
40             rbslvr,feopᵢ,rbopᵢ,xᵢ,x̂ᵢ,festatsᵢ,rbstatsᵢ
41         )
42         push!(perfs,perfᵢ)
43     end
44     return mean(perfs)
45 end
46
47 # online phase and post process
48 μₒₙ = realization(D;nparams=10,sampling=:uniform)
49 perf = local_rb_solve(rbslvr,rbop,μₒₙ)
```

FIGURE 17. ROM code used to solve the unsteady Navier Stokes application with GridapROMs

tighter tolerance of $10^{-6}$. Because shape parameters are involved, two distinct FE operators must be defined for the offline and online parameters (see lines 21 and 34). This necessitates a small adjustment to correctly execute the online phase (and post-processing), which motivates the definition of the function `local_rb_solve`. Specifically, the reduced operator `rbop` – computed at line 23 – stores the offline FE operator, which must be replaced with the appropriately defined online FE operator. This replacement is carried out at lines 34-36. In the absence of shape parameterizations, however, the standard `solve` function can be used for the online phase, which essentially performs the following instructions:

- Clustering the online parameters based on their distance to the 8 offline centroids.
- Iterating over the clustered online parameters in a *for* loop.
- Selecting a local reduced operator (see line 36) by calling `get_local` with inputs `rbop` and one representative of the cluster (e.g., the first parameter in the cluster).
- Solving the reduced equations (see line 38).

Note that this approach is valid only when the centroids used for the local subspaces and local hyper-reductions are identical, as in this test case. If this is not the case (e.g., if different numbers of centroids are selected for the two steps), then the online parameters should not be clustered and must instead be passed individually to the *for* loop. The output of the script reports an error of $2 \cdot 10^{-3}$, with speedups of $0.64$ in time and $0.48$ in memory.